\documentclass[12pt]{article}
\usepackage{amsmath,amssymb,amsthm, amsfonts}
\usepackage{etaremune, enumerate, float, verbatim}
\usepackage{algorithm}
\usepackage{algorithmic}
\usepackage[colorlinks=true, allcolors=black]{hyperref}
\usepackage{graphicx}
\usepackage{url}
\usepackage[mathlines]{lineno}
\usepackage{dsfont} 
\usepackage{tikz, graphicx, subcaption, caption} 
\usepackage[algo2e,ruled,vlined]{algorithm2e}
\usepackage{mathrsfs}
\usepackage{soul}
\usepackage{float}
\usepackage{color}
\usepackage{multirow}
\usepackage{hhline}
\definecolor{red}{rgb}{1,0,0}

\definecolor{blue}{rgb}{0,0,.9}

\definecolor{green}{rgb}{0,.6,0}

\definecolor{purp}{rgb}{.5,0,.5}

\definecolor{grey}{rgb}{.3,.3,.3}

\numberwithin{figure}{section}   
\numberwithin{table}{section}   
\numberwithin{equation}{section}   

\usepackage{tikz}
\usetikzlibrary{fit}
\usetikzlibrary{backgrounds}
\usetikzlibrary{shapes.geometric}
\usetikzlibrary{calc}
\tikzstyle{vertex}=[circle, draw=black, thick, inner sep=0pt, minimum size=6pt]
\tikzstyle{Bvertex}=[circle, black, fill, draw, inner sep=0pt, minimum size=6pt]
\tikzstyle{vtx}=[circle, white, fill, draw=black, thick, inner sep=0pt, minimum size=6pt]
\tikzstyle{gvertex}=[circle, green, fill, draw=black, inner sep=0pt, minimum size=6pt]

\newtheorem{thm}{Theorem}[section]
\newtheorem{cor}[thm]{Corollary}
\newtheorem{lem}[thm]{Lemma}
\newtheorem{prop}[thm]{Proposition}

\theoremstyle{definition}
\newtheorem{rem}[thm]{Remark}

\theoremstyle{definition}

\theoremstyle{definition}
\newtheorem{ex}[thm]{Example}

\newcommand{\Z}{\operatorname{Z}}

\newcommand{\Zp}{\operatorname{Z}_+}

\newcommand{\pd}{\gamma_P} 
\newcommand{\ptz}{\operatorname{pt}_{\Z}}
\newcommand{\ptx}{\operatorname{pt}_Y}
\newcommand{\ptp}{\operatorname{pt}_+}

\newcommand{\ppt}{\operatorname{pt}_{\pd}}
\newcommand{\capt}{\operatorname{capt}}

\newcommand{\thx}{\operatorname{th}_Y}
\newcommand{\thz}{\operatorname{th}_{\Z}}
\newcommand{\thp}{\operatorname{th}_+}

\newcommand{\thc}{\operatorname{th}_c}
\newcommand{\thpd}{\operatorname{th}_{\pd}}
\newcommand{\thxx}{\operatorname{th}_Y^\times}
\newcommand{\thxa}{\operatorname{th}_Y^\ast}
\newcommand{\thzx}{\operatorname{th}_{\Z}^\times}
\newcommand{\thza}{\operatorname{th}_{\Z}^\ast}
\newcommand{\thpx}{\operatorname{th}_+^\times}
\newcommand{\thpa}{\operatorname{th}_+^\ast}
\newcommand{\thcx}{\operatorname{th}_c^{\x}}
\newcommand{\thca}{\operatorname{th}_c^\ast}
\newcommand{\thpda}{\operatorname{th}_{\pd}^\ast}
\newcommand{\thpdx}{\operatorname{th}_{\pd}^\times}

\newcommand{\rd}{\operatorname{rd}}

\newcommand{\Gc}{\overline{G}}

\newcommand{\x}{\times}

\newcommand{\bit}{\begin{itemize}}
\newcommand{\eit}{\end{itemize}}
\newcommand{\ben}{\begin{enumerate}}
\newcommand{\een}{\end{enumerate}}
\newcommand{\beq}{\begin{equation}}
\newcommand{\eeq}{\end{equation}}
\newcommand{\bea}{\begin{eqnarray*}} 
\newcommand{\eea}{\end{eqnarray*}}
\newcommand{\bpf}{\begin{proof}}
\newcommand{\epf}{\end{proof}\ms}
\newcommand{\bmt}{\begin{bmatrix}}
\newcommand{\emt}{\end{bmatrix}}
\newcommand{\ms}{\medskip}

\newcommand{\cp}{\, \Box\,}

\newcommand{\lf}{\left\lfloor}
\newcommand{\rf}{\right\rfloor}
\newcommand{\lp}{\!\left(}
\newcommand{\rp}{\right)}

\newcommand{\du}{\,\dot{\cup}\,}
\newcommand{\noi}{\noindent}
\newcommand{\ol}{\overline}

\newcommand{\fc}{\mathsf F}
\newcommand{\ls}{\mathsf L}
\newcommand{\rn}{\mathsf R}

\title{Sharp Nordhaus-Gaddum bounds for throttling}
\author{Ryan Blair\thanks{Department of Mathematics and Statistics, California State University Long Beach, Long Beach, CA  90840, USA (ryan.blair@csulb.edu)}\and Gabriel Elvin\thanks{Department of Mathematics, California State University, San Bernardino, San Bernardino, CA 92407, USA (Gabriel.Elvin@csusb.edu)}\and Veronika Furst\thanks{Department of Mathematics, Fort Lewis College, Durango, CO  81301, USA (furst\_v@fortlewis.edu)}\and Leslie Hogben\thanks{American Institute of Mathematics, Pasadena, CA 91125, USA (hogben@aimath.org); Department of Mathematics, Iowa State University, Ames, IA 50011, USA; Department of Mathematics, Purdue University, West Lafayette, IN 47906, USA.} 
\and Tony W. H. Wong \thanks{Department of Mathematics, Kutztown University of Pennsylvania, Kutztown, PA 19530, USA (wong@kutztown.edu)}}

\begin{document}

\maketitle

\begin{abstract} 
Throttling is a graph optimization problem, where the throttling number of a graph is the minimum sum or minimum product of the number of vertices in an initial set and the time required to complete a certain graph operation. A Nordhaus-Gaddum bound refers to an upper or lower bound of the sum or product of a graph parameter together with that of its complement. In this paper, we study the Nordhaus-Gaddum sum and product bounds of the various throttling numbers (sum throttling and product throttling with or without initial cost). Graph operations considered are standard zero forcing, positive semidefinite forcing, power domination, and Cops and Robbers.
\end{abstract}

 \noi {\bf Keywords}  throttling, Nordhaus-Gaddum, propagation time, zero forcing, positive semidefinite forcing, power domination, Cops and Robbers, graph searching

\noi{\bf AMS subject classification} 05C69, 05C50, 68R10 \ms


\section{Introduction}

 A \emph{Nordhaus-Gaddum (NG) bound} is an upper or lower bound on the sum or product of the values of a parameter on a graph and its complement. Aouchiche and Hansen provide an  extensive survey of NG bounds in \cite{NGsurvey13}.   By simultaneously limiting the possible values of a graph invariant of $G$ and of its complement $\ol G$, Nordhaus-Gaddum bounds provide useful structural information about a graph, given such information about its complement.  In this paper, we focus on the graph parameters known as throttling numbers.

Throttling in graphs belongs in the intersection of graph theory and combinatorial optimization.  It involves the analysis of the efficiency of a task by minimizing the sum or the product of a measure of the resources used and the time needed to accomplish the task.  We consider two metrics:  in sum throttling, the total we aim to minimize is the sum of the amount of resources used and the time needed, while product throttling aims to minimize their product.  Small differences in the initial set of resources or in the propagation time can have a much larger overall impact on product throttling than they would on sum throttling.  Product throttling is separated into two variants, one in which initial cost is taken into account and one in which it is not but at least one time step is required; both stipulations have the effect of ruling out the trivial situation in which the product throttling number of a graph equals $0$.  

In this paper, we calculate NG sum and product bounds for sum and both types of product throttling for the following four graph parameters:  the zero forcing number $\Z(G)$, the positive semindefinite zero forcing number $\Zp(G)$, the power domination number $\pd(G)$, and the cop number $c(G)$ of a graph $G$. Tables \ref{table:summarysum} and \ref{table:summaryprod} summarize the NG sum and NG product bounds, respectively.  Over all graphs, the NG sum bounds for $\Z$ were known as was the NG sum upper bound for $\Zp$; otherwise the results are new. We provide all necessary definitions and some background in the remainder of this section.

Section \ref{s:Z-NG-th} contains our results for standard zero forcing.  We construct several novel examples (e.g., Example \ref{Rn}) used to prove that the new NG bounds we present (e.g., Theorems \ref{NG-Z-con-lower} and \ref{NG-Z-prod-lower}) are sharp, almost sharp, or tight, even in the subproblem that restricts to graphs for which both $G$ and $\ol{G}$ are connected.  In Section \ref{s:PSD-NG-th} we turn our attention to positive semidefinite zero forcing, improving known bounds and establishing new ones.  The key findings in this section are Corollaries \ref{NG-Zp-prodth-both-arb} for arbitrary graphs and \ref{NG-Zp-prodth-both} when both $G$ and $\ol{G}$ are connected, with sharpness/tightness provided by several examples.  Finally, since both power domination and Cops and Robbers are closely related to graph domination, we consider both of these processes together in Section \ref{s:PD+CR-NG-th}.  The upper and lower bounds in Theorems \ref{NG-pd-cr} and \ref{NG-pd-cr-con} comprise the main results of this section, once again relying on a series of examples to demonstrate sharpness/tightness.


\subsection{Notation for graphs and asymptotics}

A graph is simple, undirected,  finite, and has a nonempty vertex set. For a graph $G$, $V(G)$ denotes the set of vertices and $E(G)$ denotes the set of edges; $|V(G)|$ is the \emph{order} of $G$,  and an edge between vertices $u$ and $v$ is denoted by $uv$ or $vu$. 
The \emph{complement} of $G$, denoted by $\ol{G}$, is the graph $(V(G), \overline{E(G)})$ where  $uv\in \overline{E(G)}$ if and only if $uv\notin E(G)$ for $u,v\in V(G)$ and $u\ne v$. 
 Any graph whose sets of vertices and edges are subsets of $V(G)$ and $E(G)$, respectively, is a \emph{subgraph} of $G$.
Given a set $W\subseteq V(G)$,  the \emph{induced subgraph}  $G[W]$ has $V(G[W])=W$ and $E(G[W])=\{uv: uv\in E(G)\mbox{ and }u,v\in W\}$.
The \emph{union} of $G$ and $G'$ is the graph $G\,\cup\, G' = (V(G)\cup V'(G), E(G)\cup E'(G))$.   
If $ V(G)\cap V'(G)=\emptyset$, then the union is \emph{disjoint} and can be denoted by $G\du G'$. Two graphs $G$ and $G'$ are \emph{isomorphic}, denoted by $G \cong G'$, if there exists a bijection $\varphi: V(G) \to V(G')$ such that $uv\in E(G)$ if and only if $\varphi(u)\varphi(v) \in E(G')$.

Given $v\in V(G)$, the \emph{(open) neighborhood} of $v$ in $G$ is $N_G(v)=\{u: uv\in E(G)\}$, \st{and} the \emph{closed neighborhood} of $v$ in $G$ is $N_G[v]=N(v)\cup \{v\}$, and $N_G[W]=\cup_{w\in W}N_G[W]$ for $W\subseteq V(G)$; if the graph is clear from context, the subscript may be omitted.  
The \emph{degree} of a vertex $v$ is $\deg_G(v) = |N_G(v)|$.
A vertex $\ell$ of a graph $G$ is    a \emph{leaf} if $\deg_G(\ell)=1$.   A vertex $u$ of a graph $G$ is   \emph{universal} if $u$ is adjacent to every other vertex of $G$. An \emph{independent} set    is a  set of vertices of which no two  are adjacent. The \emph{independence number} of $G$, denoted by $\alpha(G)$,  is the maximum cardinality of an independent set. 

A \emph{path} in a graph $G$ is a sequence of distinct vertices $v_1,v_2,\dots,v_r$ such that for each $i$ with $1 \leq i \leq r-1$ we have $v_iv_{i+1} \in E(G)$.   
  A graph $G$ is  \emph{connected} if for each pair of vertices $u,w \in V(G)$ there exists a path $v_1,v_2,\dots,v_r$ with $v_1=u$ and $v_r=w$; a graph is \emph{disconnected} if it is not connected.   The \emph{components} of $G$ are its maximal connected subgraphs.  
  
  The \emph{path  graph} $P_n$ is the graph with $V(P_n)=\{v_1,\dots,v_n\}$ and $E(P_n)=\{{ v_iv_{i+1}}: i=1,\dots,n-1\}$. The \emph{cycle graph} $C_n$ is the graph obtained from the path graph $P_n$ just described by adding the edge $v_1v_n$. 
The \emph{complete graph} of order $n$, which has an edge between every pair of vertices,  is denoted by $K_n$.
The \emph{complete bipartite graph} $K_{p,q}$ is the graph with $V(K_{p,q})=X\du Y$ such that $|X|=p$ and $|Y|=q$ and $E(K_{p,q})=\{xy: x\in X\mbox{ and }y\in Y\}$; $K_{1,n-1}$ is called a \emph{star}. The \emph{complete multipartite graph} $K_{n_1,\dots,n_r}$ is the graph whose vertices can be partitioned into $r$ parts of sizes $n_1,\ldots,n_r$ and  the set of edges is exactly the set of  all  possible edges between distinct parts.

Let $f$ and $g$ be positive real valued functions of $n\in \mathbb{N}$ for $n\ge 2$. 
 We say $f$ is $o(g)$  if $\lim_{n\to\infty}\frac{f(n)}{g(n)}=0$.
 We say  $f$ is $O(g)$  if there exist positive constants $C,k$ such that $f(n)\le C g(n)$ for all $n\ge k$,
and  $f$ is $\Omega(g)$  if  there exist positive constants $c,k$ such that $f(n)\ge c g(n)$ for all $n\ge k$.
Finally,  $f$ is $\Theta(g)$  if $f$ is both $O(g)$ and $\Omega(g)$. 
\subsection{Parameters}\label{ss:param}

In this section we define the parameters for which we will discuss Nordhaus-Gaddum problems.
Standard zero forcing, positive semidefinite (PSD) zero forcing, and power domination are  propagation processes on a graph, where the goal is to fill  all the vertices  (starting with  each vertex filled or unfilled); unfilled vertices are filled by applying a {color change rule}.\footnote{Originally colors were used for the vertices, leading to the name `color change rule', but filled/unfilled has  become standard due to its suitability for all media.} 
Standard and PSD zero forcing originated in combinatorial matrix theory, providing upper bounds on the nullities of certain symmetric matrices whose off-diagonal pattern of nonzero entries is described by the given graph; 
standard zero forcing was also introduced in control of quantum systems. Power domination originated from the problem of optimal placement of PMUs, monitoring devices in an electrical network. More information on the origins of these parameters and relevant references can be found in \cite{HLSbook}.

The \emph{standard color change rule} is: If $w$ is the unique unfilled neighbor of a filled vertex $v$, then  fill  $w$.  
The  \emph{PSD color change rule} is: Let $B$ be the set of currently filled vertices and let $W_1,\dots, W_k$ be the sets of vertices of the  components of $G[V(G)\setminus B]$.  If $v\in B$, $w\in W_i$, and $w$ is the only unfilled neighbor of  $v$ in $G[W_i\cup B]$, then fill $w$. 
It is possible that there is only one component of $G[V(G)\setminus B]$, and in that case the effect of the PSD color change rule is the same as that of the standard color change rule.   
Forcing using the standard color change rule or the PSD color change rule is also called \emph{standard zero forcing} or \emph{PSD forcing}.
 Repeated application of the standard or PSD color change rule until no more vertices can be filled is called the \emph{standard zero forcing propagation process} or the \emph{PSD forcing propagation process}.
\emph{Power domination} begins with a \emph{domination step}, in which every neighbor of a filled vertex is filled.  After the domination step, the standard color change rule is applied; together these steps are  called the \emph{power domination propagation process}.  

 A set $B$ of vertices is called a \emph{standard zero forcing set, PSD forcing set}, or \emph{power dominating set}, respectively, if starting with the vertices in $B$ filled and the remaining vertices unfilled, the respective propagation process can fill all the vertices.  The \emph{standard zero forcing number} $\Z(G)$ is the minimum cardinality of a standard zero forcing set, and similarly for the \emph{PSD forcing number} $\Zp(G)$ and \emph{power domination number} $\pd(G)$.
 Power domination has natural connections to domination, so we provide definitions and notation for that also: A set $S\subseteq V(G)$ is a \emph{dominating set of $G$} if every vertex of $G$ is in $S$ or is a neighbor of a vertex in $S$; the \emph{domination number}  
 $\gamma(G)$ is the minimum cardinality of a dominating set.  

 \emph{Cops and Robbers} is a  pursuit-evasion game played on a graph by two players. One player places and moves a set of cops and the other places and moves one robber. The goal for the cops is to catch the robber by occupying the same vertex the robber occupies, which is called \emph{capture}.
The goal of the robber is to avoid capture. Initially the cops are placed  on a multiset  $B$ of vertices (more than one cop can occupy a single vertex), and then the robber is placed on a vertex. For the rest of the game, players alternate turns, starting with the cops player; a \textit{round} consists of one turn for the cops player and one for the robber player. A turn consists of repositioning all cops or the robber so that each is at most distance one away from its current vertex. The cops win the game if a cop captures the robber after a finite number of rounds, whereas if the robber has a strategy to evade the cops indefinitely, the robber wins. If choosing the set $B$ results in  the cops player  being able to win   for any choice of  strategy by the robber player, then $B$ is called a \emph{capture set}. The \emph{cop number} $c(G)$ of a graph $G$ is the minimum cardinality of a capture set; see \cite{CRbook} for more information on Cops and Robbers.

We often use $Y$ to represent one or more of  the processes such as standard zero forcing,  PSD forcing, power domination, or Cops and Robbers,  or the $Y$-number  $Y(G)$ of a graph $G$, a $Y$-set, etc.

We define $Y$-propagation time for $Y$ either standard zero forcing  or  PSD forcing as follows: 
Start with $B^{[0]}=B^{(0)}=B\subseteq V(G)$.
Assume $B^{(i)}$ and $B^{[i]}$ have been constructed. Then 
\[B^{(i+1)}=\{w: 
\mbox{ $w$ can be  filled (given that all vertices in $B^{[i]}$  are filled)}\} .\] 
and $B^{[i+1]}=B^{[i]}\cup B^{(i+1)}$. 
 For power domination, $B^{[1]}=N[B]$  and $B^{(1)}=B^{[1]}\setminus B$, and   $B^{(i+1)}$ and $B^{[i+1]}$ are defined as above for $i\ge 1$.
The   \emph{$Y$-propagation time}  (for standard or PSD forcing or power domination) of $B\subseteq V(G)$, denoted by $\ptx(G;B)$, is the least $t$ such that $B^{[t]}=V(G)$; if $B^{[t]}\ne V(G)$ for all $t$, then $\ptx(G;B)=\infty$. 
 The \emph{$k$-propagation time of $G$} for $Y$ is
$ \ptx(G,k)=\min_{|B|= k}\ptx(G;B).$  The \emph{$Y$-propagation time of  $G$}  is 
$\ptx(G)=\ptx(G,Y(G)).$

For Cops and Robbers, there is no analog of $\ptx(G;B)$ because the cops can be placed anywhere.  The  \emph{$k$-capture time} {$\capt_k(G)$} is the minimum number of rounds needed for $k\ge c(G)$ cops to capture the robber on $G$ over all possible games (when both the cops and the robber play optimally); $\capt_k(G)=\infty$ for $k< c(G)$.  
The \emph{capture time} of $G$ is the $c(G)$-capture time and is denoted by {$\capt(G)$}.  

 For $Y$ one of standard zero forcing,  PSD forcing, or power domination, we define the sum and two types of product throttling number of a set, the $k$-throttling number, and the throttling number of the graph.  Let $G$ be a graph of order $n$. 

\emph{Sum throttling\footnote{This is often just called  the throttling number in the literature.} for $Y$}: For a set $B\subseteq V(G)$ and $1\le k\le n$, 
\[\thx(G;B)=|B|+\ptx(G;B), \  \thx(G,k)=k+\ptx(G,k), 
\mbox{ 
and }\thx(G)=\min_{1\le k\le n}\thx(G,k).\]

\emph{No initial cost product throttling for $Y$:}   It is assumed  that $G$  has an edge.  For a set $B\subsetneq V(G)$ and $1\le k\le n-1$, 
\[\thxa(G;B)=|B|\ptx(G;B), \  \thxa(G,k)=k \cdot\ptx(G,k), 
\mbox{ 
and }\thxa(G)=\min_{1\le k\le n-1}\thxa(G,k).\]

\emph{Initial cost product throttling for $Y$:}   For   a set $B\subseteq V(G)$ and $1\le k\le n$, 
\[\thxx(G;B)=|B|(1+\ptx(G; B)), \ 
\thxx(G,k)=k (1+\ptx(G,k)), 
\mbox{ 
and }\thxx(G)=\min_{1\le k\le n}\thxx(G,k).\]

For Cops and Robbers we define the $k$-throttling number  of the graph by replacing $\ptx(G,k)$ by $\capt_k(G)$ and then using the same definitions for throttling.


\subsection{Nordhaus-Gaddum bounds for sum and product throttling}\label{s:NG}

In this section, we present some elementary observations and summarize the results of this paper.  Table~\ref{table:summarysum} lists Nordhaus-Gaddum sum bounds for the graph searching  parameters and throttling parameters studied here.  Note that when the graphs considered are arbitrary, the NG sum bounds for standard zero forcing sum throttling and the NG sum upper bound for PSD forcing sum throttling were previously known in the literature, as was a weaker NG sum lower bound for PSD  forcing sum throttling; all the remaining bounds are new.  Table~\ref{table:summaryprod} presents new NG product bounds for all forcing parameters and throttling parameters under consideration.  The following remark contains a useful collection of some elementary observations.

\begin{rem}\label{NG-u} \label{NG-u2} Let $G$ be a graph (and assume each of $G$ and $\ol G$ has an edge for no initial cost product throttling), and let $Y$ be one of the parameters $\Z, \Zp, \pd$, or $c$.  Then $1\le Y(G)\le n$, and if $G$ has an edge, then $Y(G)\le n-1$.  
If $G$ does not have an edge, then $\thx(G)=\thxx(G)= n$ (and $\thxa(G)$ is not defined).

If $G$ has an edge, then $\thx(G),\thxx(G)\le n$ and $\thxa(G)\le n-1$;
furthermore, $\thx(G)\ge Y(G)+1  \ge 2$, $\thxa(G)\ge Y(G)  \ge 1$, and  $\thxx(G)\ge \min\{2Y(G), n\}  \ge 2$. 
As noted in \cite{product2}, $\thxx(G)>\thxa(G)\ge \thx(G)-1$, which also implies  $\thxx(G)\ge \thx(G)$.  
Thus 
\[4\le \thx(G)+\thx(\ol G)\le 2n, \ 2\le \thxa(G)+\thxa(\ol G)\le 2n-2, \mbox{ and } 4\le \thxx(G)+\thxx(\ol G)\le 2n.\] 
\[4\le \thx(G)\cdot\thx(\ol G)\le n^2, \ 1\le \thxa(G) \cdot \thxa(\ol G)\le (n-1)^2, \mbox{ and } 4\le \thxx(G) \cdot \thxx(\ol G)\le n^2.\] 
\end{rem}

Given a graph parameter $\zeta$, a  \emph{sharp Nordhaus-Gaddum sum bound} is a bound on $\zeta(G)+\zeta(\Gc)$ that is attained by infinitely many graphs.  Analogously, a \emph{sharp Nordhaus-Gaddum product bound} is a bound on $\zeta(G)\cdot\zeta(\Gc)$ that is attained by infinitely many graphs.  We say a Nordhaus-Gaddum sum (respectively, product)  bound $b(n)$ on $\zeta$ is \emph{almost sharp} if   infinitely many graphs satisfy $\zeta(G)+\zeta(\Gc)=b(n)\pm k$ (respectively, $\zeta(G)\cdot \zeta(\Gc)=b(n)\pm k$) for some  positive integer $k< 10$.   A Nordhaus-Gaddum sum (respectively, product) bound $b(n)\pm o(b(n))$ is \emph{tight} if infinitely many graphs satisfy $\zeta(G)+\zeta(\Gc)=b(n)\pm o(b(n))$ (respectively, $\zeta(G)\cdot \zeta(\Gc)=b(n)\pm o(b(n))$). 

The NG problem with the restriction  that both $G$ and $\ol G$ must be connected has also been studied for many parameters (see \cite{PD-NG-REUF} for work on power domination and \cite{NGsurvey13} for various parameters).  Since throttling is usually studied for connected graphs (and the upper bound is often trivially realized using isolated vertices), we emphasize this restricted version.
 
\begin{table}[H]
\begin{center}
\renewcommand{\arraystretch}{1.2}
\footnotesize\begin{tabular}{|c|c|c|c|c|c|c|}
\hline
\multirow{2}{*}{$Y$}& \multicolumn{2}{c|}{NG sum bounds on $\thx$}& \multicolumn{2}{c|}{NG sum bounds on $\thxa$}& \multicolumn{2}{c|}{NG sum bounds on $\thxx$}\\
\cline{2-7}
& Lower& Upper& Lower& Upper& Lower& Upper\\
\hhline{|=|=|=|=|=|=|=|}

\multirow{4}{*}{$\Z$}& \multirow{2}{*}{(C,A) $n+o(n)$}& (C) $2n-2$& & (C) $2n-4$& \multicolumn{2}{c|}{\multirow{4}{*}{N/A}}\\
& \multirow{2}{*}{Rmk \ref{r:thz-LBsum}}& Prop~\ref{NG-Z-con}& (C,A) $\frac{5n}{4}-2$\; AS & Prop~\ref{NG-Z-con}& \multicolumn{2}{c|}{}\\
\cline{3-3}\cline{5-5}
& \multirow{2}{*}{Thm~\ref{t:NG-sum:Z}}& (A) $2n$& Thm~\ref{NG-Z-con-lower}& (A) $2n-2$& \multicolumn{2}{c|}{}\\
& & Thm~\ref{t:NG-sum:Z}& & Prop~\ref{NG-Z}& \multicolumn{2}{c|}{}\\
\hline

\multirow{4}{*}{$\Zp$}& \multirow{2}{*}{(C,A) $n$\; AS} & (C) $2n-o(n)$& \multirow{2}{*}{(C,A) $n-2$\; AS}& (C) $2n-o(n)$& \multirow{2}{*}{(C,A) $n$\; AS}& \multirow{2}{*}{(C,A) $2n$}\\
& \multirow{2}{*}{Cor~\ref{NG-Zp-prodth-both}}& Cor~\ref{NG-Zp-prodth-both}& \multirow{2}{*}{Cor~\ref{NG-Zp-prodth-both}}& Cor~\ref{NG-Zp-prodth-both}& \multirow{2}{*}{Cor~\ref{NG-Zp-prodth-both}}& \multirow{2}{*}{Cor~\ref{NG-Zp-prodth-both}}\\
\cline{3-3}\cline{5-5}
& \multirow{2}{*}{Cor~\ref{NG-Zp-prodth-both-arb}}&  (A) $2n$& \multirow{2}{*}{Cor~\ref{NG-Zp-prodth-both-arb}}& (A) $2n-2$& \multirow{2}{*}{Cor~\ref{NG-Zp-prodth-both-arb}}& \multirow{2}{*}{Cor~\ref{NG-Zp-prodth-both-arb}}\\
& & Cor~\ref{NG-Zp-prodth-both-arb}& & Cor~\ref{NG-Zp-prodth-both-arb}& &\\
\hline

\multirow{4}{*}{$\pd$}& (C) $6$& (C)$^\dagger$ $\lf\frac{n}{3}\rf+5$& (C) $4$& (C)$\lf\frac{n}{2}\rf+2$& \multirow{2}{*}{(C,A) $6$}& (C)$^\ddagger$ $\frac{6n}{7}+3$\\
& Thm~\ref{NG-pd-cr-con}& Thm~\ref{NG-pd-cr-con}& Thm~\ref{NG-pd-cr-con}& Thm~\ref{NG-pd-cr-con}& \multirow{2}{*}{Thm~\ref{NG-pd-cr-con}}& Thm~\ref{NG-pd-cr-con}\\
\cline{2-5}\cline{7-7}
& (A) $5$& (A) $n+2$& (A) $3$& (A) $n$& \multirow{2}{*}{Thm~\ref{NG-pd-cr}}& (A) $n+3$ \\ 
& Thm~\ref{NG-pd-cr}& Thm~\ref{NG-pd-cr}& Thm~\ref{NG-pd-cr}& Thm~\ref{NG-pd-cr}& &Thm~\ref{NG-pd-cr}\\
\hline

\multirow{5}{*}{$c$}& \multirow{2}{*}{(C) $6$}& (C) $\min\Big(\lf\frac{n}{4}\rf+4,$& \multirow{2}{*}{(C) $4$}& \multirow{2}{*}{(C) $\lf\frac n 2\rf +2$ ?}& & \multirow{2}{*}{(C) $n+4$ ?}\\
& \multirow{2}{*}{Thm~\ref{NG-pd-cr-con}}& $\frac{2n}{(\log n)^{\frac{1}{2}-o(1)}}\Big)$ ?& \multirow{2}{*}{Thm~\ref{NG-pd-cr-con}}& \multirow{2}{*}{Thm~\ref{NG-pd-cr-con}}& (C,A) $6$& \multirow{2}{*}{Thm~\ref{NG-pd-cr-con}}\\
& &Thm~\ref{NG-pd-cr-con}& & & Thm~\ref{NG-pd-cr-con}&\\
\cline{2-5}\cline{7-7}
& (A) $5$& (A) $n+2$& (A) $3$& (A) $n$& Thm~\ref{NG-pd-cr}&(A) $n+4$\\
& Thm~\ref{NG-pd-cr}& Thm~\ref{NG-pd-cr}&Thm~\ref{NG-pd-cr}&Thm~\ref{NG-pd-cr}& &Thm~\ref{NG-pd-cr}\\
\hline
\end{tabular}
\end{center}
\caption{
Summary of the Nordhaus-Gaddum sum bounds for throttling  for graphs of order $n\ge 5$ except as indicated by  $\dagger$  ($n\ge 12$) or $\ddagger$ ($n\ge 16$). More precise information is sometimes given in the cited result.\\  
The label (C) indicates that both $G$ and $\ol{G}$ are connected.\\
The label (A) indicates that both $G$ and $\ol{G}$ are arbitrary graphs.\\
The label (C,A) indicates that the bound holds for both types of graph; when two references are listed, the first  is for connected graphs.\\ 
All asymptotic bounds are tight.  Exact bounds marked with ``?'' may not even be tight,  exact bounds marked with  ``AS" are almost sharp, and  exact bounds with no additional markings  are sharp.
}
\label{table:summarysum}
\end{table}


\begin{table}[H]
\begin{center}
\renewcommand{\arraystretch}{1.2}
\footnotesize\begin{tabular}{|c|c|c|c|c|c|c|}
\hline
\multirow{2}{*}{$Y$}& \multicolumn{2}{c|}{NG product bounds on $\thx$}& \multicolumn{2}{c|}{NG product bounds on $\thxa$}& \multicolumn{2}{c|}{NG product bounds on $\thxx$}\\
\cline{2-7}
& Lower& Upper& Lower& Upper& Lower& Upper\\
\hhline{|=|=|=|=|=|=|=|}

\multirow{4}{*}{$\Z$}& & (C) $(n-1)^2$& & (C) $(n-2)^2$& \multicolumn{2}{c|}{\multirow{4}{*}{N/A}}\\
& (C,A) $\Theta(n^{3/2})$ 
& Prop~\ref{NG-Z-con}
& (C,A) $\frac{3n^2}{8}-O(n)$& Prop~\ref{NG-Z-con}& \multicolumn{2}{c|}{}\\
\cline{3-3}\cline{5-5}
& Cor \ref{cor:theta} 
& (A) $n^2$& Thm~\ref{NG-Z-prod-lower}& (A) $(n-1)^2$& \multicolumn{2}{c|}{}\\
& & Prop~\ref{NG-Z}& & Prop~\ref{NG-Z}& \multicolumn{2}{c|}{}\\
\hline

\multirow{4}{*}{$\Zp$}&  (C) $ 3n-9$\; AS & (C) $n^2-o({n^2})$&  (C) $ 2n-8$\; AS& (C) $n^2-o({n^2})$&  (C) $3n$& \multirow{2}{*}{(C,A) $n^2$}\\
& Cor~\ref{NG-Zp-prodth-both} &  Cor~\ref{NG-Zp-prodth-both}& Cor~\ref{NG-Zp-prodth-both}&  Cor~\ref{NG-Zp-prodth-both}& Cor~\ref{NG-Zp-prodth-both} & \multirow{2}{*}{Cor~\ref{NG-Zp-prodth-both}}\\
\cline{2-6}
& (A) $2n$&  (A) $n^2$& (A) $n-1$& (A) $(n-1)^2$& (A) $ 2n$& \multirow{2}{*}{Cor~\ref{NG-Zp-prodth-both-arb}}\\
& Cor~\ref{NG-Zp-prodth-both-arb}& Cor~\ref{NG-Zp-prodth-both-arb}& Cor~\ref{NG-Zp-prodth-both-arb}& Cor~\ref{NG-Zp-prodth-both-arb}& Cor~\ref{NG-Zp-prodth-both-arb}&\\
\hline

\multirow{4}{*}{$\pd$}& (C) $9$& (C) $n+o(n)$& (C) $4$& \multirow{2}{*}{(C,A) $n$}& (C) $9$& (C)  $\Theta(n)$\\
& Thm~\ref{NG-pd-cr-con}& Thm~\ref{NG-pd-cr-con}& Thm~\ref{NG-pd-cr-con}& \multirow{2}{*}{Thm~\ref{NG-pd-cr-con}}& Thm~\ref{NG-pd-cr-con}& Cor~\ref{NG-pd-theta-con}\\
\cline{2-4}\cline{6-7}
& (A) $6$& (A) $2n$& (A) $2$& \multirow{2}{*}{Thm~\ref{NG-pd-cr}}& (A) $8$& (A) $ \Theta(n)$ \\
& Thm~\ref{NG-pd-cr}& Thm~\ref{NG-pd-cr}& Thm~\ref{NG-pd-cr}& & Thm~\ref{NG-pd-cr}&{Cor~\ref{NG-pd-cr-theta}}\\
\hline

\multirow{4}{*}{$c$}& (C) $9$& (C) {$n+\lf\frac{n}{2}\rf+3$} ?& (C) $4$& (C) { $n$} ?& (C) $9$& (C)  $ 4n$ ?\\
& Thm~\ref{NG-pd-cr-con}& Thm~\ref{NG-pd-cr-con}& Thm~\ref{NG-pd-cr-con}& { Thm~\ref{NG-pd-cr-con}}& Thm~\ref{NG-pd-cr-con}& { Thm~\ref{NG-pd-cr-con}}\\
\cline{2-7}
& (A) $6$& (A) $2n$& (A) $2$& { (A) $n$}& (A) $8$& { (A) $4n$}\\
& Thm~\ref{NG-pd-cr}& Thm~\ref{NG-pd-cr}& Thm~\ref{NG-pd-cr}& { Thm~\ref{NG-pd-cr}}& Thm~\ref{NG-pd-cr}& { Thm~\ref{NG-pd-cr}}\\
\hline
\end{tabular}
\end{center}
\caption{Summary of the Nordhaus-Gaddum product bounds for throttling  for graphs of order $n\ge6$.  More precise information is sometimes given in the cited result.\\
(C) indicates that both $G$ and $\ol{G}$ are connected.\\
(A) indicates that both $G$ and $\ol{G}$ are arbitrary graphs.\\
The label (C,A) indicates that the bound holds for both types of graph.\\ 
All asymptotic bounds are tight. 
 Exact bounds marked with ``?" may not even be tight,  exact bounds marked with ``AS" are almost sharp, and  exact bounds with no additional markings  are sharp.} 
\label{table:summaryprod}
\end{table}



\section{Standard zero forcing} \label{s:Z-NG-th}

We begin our study of Nordhaus-Gaddum problems for throttling with the most restrictive of our four forcing parameters, namely standard zero forcing.  Recall that initial cost product throttling for standard zero forcing is not studied because $\thzx(G)=n$ for every graph of order $n$ \cite{product2}. By considering the other two types of throttling, we establish the first row of Tables \ref{table:summarysum} and \ref{table:summaryprod}. The following result is known for sum throttling. 
\begin{thm}\label{t:NG-sum:Z}{\rm \cite[Theorem 11.24]{HLSbook}}
Let $G$ be a graph of order $n$.  
Then  
\[n+o(n)\le \thz(G)+\thz(\Gc)\le 2n.\]
Furthermore, the upper bound is sharp and the lower bound is tight regardless of whether both $G$ and $\Gc$ are required to have an edge.
\end{thm}

Next, we state a useful tool for studying no initial cost product throttling. Define $k(G,1)=\min\{k:\ptz(G,k)=1\}$.  
\begin{thm}\label{t:NG-NCproductZset}{\rm \cite{product2}, \cite[Theorem 11.73]{HLSbook}} For any graph $G$ of order $n\geq 2$, $\thza(G)$ is the least $k$ such
that $\ptz(G, k) = 1$, i.e., $\thza(G) = k(G, 1)$. Furthermore, $k(G, 1) \geq \frac{n}
{2}$.
\end{thm}

Propositions \ref{NG-Z} and \ref{NG-Z-con} give straightforward NG sum and NG product upper bounds.
\begin{prop}\label{NG-Z}
For any graph $G$ of order $n\geq 2$,
\ben[$(1)$]
\item $ \thz(G)\cdot \thz(\ol G) \le n^2. $
\item $\thza(G)+ \thza(\ol G) \le 2(n-1) $ and $ \thza(G)\cdot \thza(\ol G) \le (n-1)^2. $
\een
All these bounds are  sharp.
\end{prop}
\bpf The bounds follow from Remark \ref{NG-u}.  The bound for $\thz$ is sharp using $G=K_n$.
The bounds for $\thza$ are sharp using $G=K_{1,n-1}$: To see $\thza(K_{1,n-1})=n-1$, note that $\Z(K_{1,n-1})=n-2$ and $\ptz(K_{1,n-1})=2$, so $\thza(K_{1,n-1})>n-2$. 
Since $\ol{K_{1,n-1}}=K_1\du K_{n-1}$ and $\Z(K_1\du K_{n-1})=n-1$, 
we have $\thza(K_1\du K_{n-1})=n-1$.  
\epf

Next we consider the case where both $G$ and $\ol G$ are required to be connected (and the order is at least four).
\begin{rem}\label{r:thz-LBsum}
    The path is the example used in \cite{HLSbook} to prove the tightness of the lower bound  in Theorem~\ref{t:NG-sum:Z}, so that lower bound remains valid if both graphs are required to be connected.  
\end{rem}

When we require both $G$ and $\Gc$ to be connected, we get slightly different results for NG  upper bounds.  The example given in \cite{HLSbook} for the sharpness of the upper bound in Theorem~\ref{t:NG-sum:Z} is a  specific kind of cograph, which implies one of $G$ and $\Gc$ is disconnected, 
as we explain next.

A graph is a \emph{cograph}  if it does not have an induced $P_4$.  Any cograph can be constructed as a sequence of join or  disjoint union operations on sets of  isolated  vertices \cite[p.\ 266]{HLSbook}.  
If the last operation is a disjoint union, then $G$ is disconnected, and if the last operation is a join, then $\ol G$ is disconnected.  Thus when both $G$ and $\ol G$ are required to be connected, neither $G$ nor $\ol G$ can be a cograph.

\begin{prop}\label{NG-Z-con}
Let $G$ be a graph of order $n\ge 4$ such that both $G$ and $\ol G$ are connected.  Then
\ben[$(1)$]
\item\label{c:thz-lower} $ \thz(G)+ \thz(\ol G) \le 2(n-1) $ and  $ \thz(G)\cdot \thz(\ol G) \le (n-1)^2 $. 
\item $ \thza(G)+ \thza(\ol G) \le 2(n-2) $ and $ \thza(G)\cdot \thza(\ol G) \le (n-2)^2 $. \een All bounds are sharp.
\end{prop}
\bpf    
Since $G$ and $\ol G$ are both connected, neither is a cograph, and hence each must contain an induced $P_4$.  By \cite[Theorem 11.22]{HLSbook}, $\thz(G), \thz(\ol G) < n$, which gives the  bounds in \eqref{c:thz-lower}.  If $(x,y,z,w)$ is an induced path in a connected graph $H$, then $\thza(H)\le n-2$ by choosing $B=V(H)\setminus\{y,z\}$ as the initial set. Thus $\thza(G), \thza(\ol G)\le n-2$.
Example \ref{ex:star+leaf} below shows that the bounds are sharp.
\epf

\begin{figure}[h!]
     \centering
     \begin{tikzpicture}[scale=1.3]
     \draw (0:1) -- (0:0);
     \draw (360/7:1) -- (0:0);
     \draw (360/7*2:1) -- (0:0);
     \draw (360/7*3:1) -- (0:0);
     \draw (360/7*4:1) -- (0:0);
     \draw (360/7*5:1) -- (0:0);
     \draw (360/7*6:1) -- (0:0);
     \foreach \d in {0,360/7,360/7*2,360/7*3,360/7*4,360/7*5,360/7*6}{
        \foreach \r in {0,...,1}{
            \draw[fill=white] (\d:\r) circle (0.1);
        }
    }
     \draw (360/7:1.1) -- (360/7:2);
     \foreach \d in {360/7}{
        \foreach \r in {2}{
            \draw[fill=white] (\d:\r)
            circle (0.1);
        }
     }
     \draw (1.6,1.6) node{$u$};
     \draw (1,.85) node{$w$};
     \draw (0.05,-0.1) node[anchor=north]{$c$};
     \end{tikzpicture}
     \qquad
      \begin{tikzpicture}[scale=1.3]
     \draw (180:1) -- (0:1);
     \draw (-135:1) -- (45:1);
     \draw (-90:1) -- (90:1);
     \draw (-45:1) -- (135:1);
     \draw (45:1) -- (90:1) -- (135:1) -- (180:1) -- (225:1) -- (270:1) -- (315:1) -- (0:1);
     \draw (0:1) -- (90:1) -- (180:1) -- (270:1) -- (0:1);\
     \draw (45:1) -- (135:1) -- (225:1) -- (315:1) -- (45:1);
     \draw (0:1) -- (135:1) -- (270:1) -- (45:1) -- (180:1) -- (315:1) -- (90:1) -- (225:1) -- (0:1);
     \draw (0:1) -- (0:2);
     \foreach \d in {0,45,...,315}{
            \draw[fill=white] (\d:1) circle (0.1);
    }
     \draw[fill=white] (0:2) circle (0.1);
     \draw (1.2,0.2) node{$u$};
     \draw (1.05,.7) node{$w$};
     \draw (2.2,0.2) node{$c$};
     \end{tikzpicture}
\caption{The graphs $\ls_n$ and $\overline{\ls_n}$ from Example \ref{ex:star+leaf} for $n=9$. }
\label{f:NG-spider}
\end{figure}
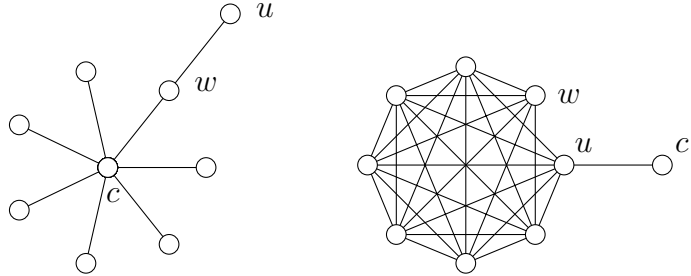

 \begin{ex}\label{ex:star+leaf} Let $\ls_n$ be the graph obtained from $K_{1,n-2}$ by adding a new leaf to one leaf of $K_{1,n-2}$; this is the spider $S(2,\underbrace{1,\dots,1}_{n-3 \ \text{times}})$, which has one leg of length 2 and $n-3$ legs of length 1.  Then $\ol{\ls_n}$ can be obtained from $K_{n-1}$ by removing an edge $e=uw$ and adding a leaf adjacent to $u$.  Both graphs are shown in Figure \ref{f:NG-spider}.  
 We have $\Z(\ls_n)=\Z(\ol{\ls_n})=n-3$, $\ptz(\ls_n)=2$, and $\ptz(\ol{\ls_n})=3$. Since $\ptz(\ls_n,n-2)=\ptz(\ol{\ls_n},n-2)=1$, we see that $\thz(\ls_n)=\thz(\ol{\ls_n})=n-1$ and  $\thza(\ls_n)=\thza(\ol{\ls_n})=n-2$.
\end{ex}

 Next we turn our attention to NG  lower bounds. Recall that Theorem \ref{t:NG-sum:Z}  and Remark \ref{r:thz-LBsum} established a tight NG sum lower bound (for both arbitrary and connected graphs).

In the proof of Theorem \ref{NGprod-lower-thz} we make use of the following definitions. Consider a fixed forcing process for a graph $G$. For each $v \in V(G)$, define the {\em round  function of $v$}, $\rd(v)$, to be  the round in which  vertex $v$ is filled; note that if $v$ is filled initially, then $\rd(v)=0$.  The set of vertices of $G$ that are filled by the $i$th round are partitioned into {\em active vertices} that have not performed a force by the $i$th round and {\em inactive vertices} that have  performed a force by the $i$th round (and thus are no longer capable of forcing).

\begin{thm}\label{NGprod-lower-thz}
    Let $G$ be a graph of order $n\ge 4$.  Then $ (2-\sqrt{3})n^{\frac{3}{2}}-(1-\sqrt{3}/2)n < \thz(G)\cdot\thz(\ol G)$.
\end{thm}

\bpf
Choose $A \subseteq V(G)$ such that $\thz(G)=\thz(G;A)$ and fix a forcing process for $G$ using $A$ as the forcing set.
Define $B^{[i]}=\{v\in V(G)|\rd(v)\leq i\}$, 
and note that $|B^{[i]}\setminus B^{[i-1]}|\leq |A|$ for every $i\geq 1$ since at most $|A|$ vertices can be forced at the $i$th round of propagation. Additionally, define $F_i \subseteq B^{[i-1]}$ to be the set of inactive vertices after the $i$th round of propagation (so $B^{[i]}\setminus F_i$ is the set of active vertices after the $i$th round). 

 Suppose first that   $|A| \geq (1 - \sqrt{3}/2)n$. Then $\thz(G) \geq (1 - \sqrt{3}/2)n + 1>(1 - \sqrt{3}/2)n$. Additionally, it is known that $\thz(H)\ge 2\sqrt n - 1$ for every graph $H$ of order $n$ \cite{BY13}, \cite[Theorem 11.14]{HLSbook}. Hence, 
{\begin{equation*}
\thz(G)\cdot\thz(\ol G)> (1-\sqrt{3}/2)n \cdot(2\sqrt n - 1) =(2-\sqrt{3})n^{\frac{3}{2}}-(1-\sqrt{3}/2)n.
\end{equation*}}

 Now suppose that $|A| < (1 - \sqrt{3}/2)n$. Let $m$ be the smallest integer such that $|B^{[m]}|\geq \frac{n}{2}$. Since $|B^{[m-1]}|< \frac{n}{2}$ and $|B^{[m]}\setminus B^{[m-1]}|\leq |A|$, we have that  $|B^{[m]}|<\frac{n}{2}+(1-\sqrt 3 / 2)n = \left(\frac{3  - \sqrt 3}{2}\right)n$  and $|\ol{B^{[m]}}|> \left(\frac{\sqrt 3 - 1}{2}\right)n$. It is known that the number of active vertices each round is constant and equal to the number of vertices in the initial forcing set  \cite[p. 188]{HLSbook}. Hence, $|B^{[m]}\setminus F_m|=|A|$ 
and  $|F_m|> \left(\frac{\sqrt 3 - 1}{2}\right)n$. 
Since each vertex in $F_m$ applied a force at the $m$th or earlier round of propagation to fill a vertex in $ B^{[m]}$, there are no edges in $G$ between  $F_m$ and $V(G)\setminus B^{[m]}$. Consequently, $\ol G$ contains more than $\left(\frac{\sqrt 3 - 1}{2}\right)^2 n^2 = (1 - \sqrt 3 / 2)n^2$  edges. By  \cite{BY13}, $\Z(H)>\frac{|E(H)|}{n}$ for every graph $H$ of order $n$. Hence, {$\Z(\ol G)> (1 - \sqrt 3 / 2)n$ and $\thz(\ol G)> (1 - \sqrt 3 / 2)n$}. Since we know $\thz(G)\ge 2\sqrt n - 1$, 
\[
\thz(G)\cdot\thz(\ol G)> (2\sqrt n - 1)(1 - \sqrt 3 / 2)n
 = (2-\sqrt{3})n^{\frac{3}{2}}-(1-\sqrt{3}/2)n. \qedhere\]
\epf

 \begin{cor}\label{cor:theta}
     The NG product lower bound for $\thz$ is $\Theta(n^{3/2})$ (both for arbitrary graphs and when both $G$ and $\Gc$ must be connected).
\end{cor}
\bpf
    When $n$ is a perfect square, $\thz(P_n)=2\sqrt n -1$ and $\thz(\ol{P_n})=n-1$ \cite[Theorem 11.23]{HLSbook}, so $\thz(P_n)\cdot\thz(\ol{P_n})=2n\sqrt n -n-2\sqrt n+1$. Hence, the NG product lower bound is at most $2n\sqrt n-O(n)$. The result then follows from Theorem \ref{NGprod-lower-thz}. \epf

 Next we consider   lower bounds for the NG sum and NG product  for $\thza$. The bounds  that we establish in Theorems~\ref{NG-Z-con-lower} (NG Sum) and \ref{NG-Z-prod-lower} (NG product)  are almost sharp and tight, respectively,  regardless of whether $G$ and $\Gc$ are required to be connected, and we begin with an example.

\begin{figure}[!h]
    \centering
    \begin{subfigure}[b]{0.45\textwidth}
    \centering
    \begin{tikzpicture}
    \draw[gray!50, fill=gray!50] (-2,3.5) -- (-2,2.5) -- (2,2.5) -- (2,3.5) -- cycle;
    \foreach \y in {-2,...,2}{
        \draw[white] (-2,3+0.167*\y) -- (2,3+0.167*\y);
    }
    \foreach \t in {0,4}{
        \foreach \x in {-2,...,2}{
            \draw (-2+0.167*\x+\t,3) -- (-2+0.167*\x+\t,0);
        }
    }
    \draw[fill=gray!50] (0,0) node{$K_{2m}$} ellipse (3 and 1);
    \draw[fill=white] (-2,3) node{$\ol{K_m}$} circle (1);
    \draw[fill=white] (2,3) node{$\ol{K_m}$} circle (1);
    \end{tikzpicture}
    \caption{$\rn_{4m}$}
    \label{fig:x}
    \end{subfigure}
    \begin{subfigure}[b]{0.45\textwidth}
    \centering
    \begin{tikzpicture}
    \draw[gray!50, fill=gray!50] (-2.5,2.5) -- (-1.5,2.5) -- (3,0) -- (-3,0) -- cycle;
    \draw[gray!50, fill=gray!50] (2.5,2.5) -- (1.5,2.5) -- (-3,0) -- (3,0) -- cycle;
    \foreach \y in {-2,...,2}{
        \draw (-2,3+0.167*\y) -- (2,3+0.167*\y);
    }
    \foreach \t in {0,4}{
        \foreach \x in {-2,...,2}{
            \draw[white] (-2+0.167*\x+\t,3) -- (-2+0.167*\x+\t,0);
        }
    }
    \draw[fill=white] (0,0) node{$\ol{K_{2m}}$} ellipse (3 and 1);
    \draw[fill=gray!50] (-2,3) node{$K_m$} circle (1);
    \draw[fill=gray!50] (2,3) node{$K_m$} circle (1);
    \end{tikzpicture}
    \caption{$\ol{\rn_{4m}}$}
    \label{fig:enter-label}
    \end{subfigure}
    \caption{ The graph $\rn_{4m}$ and its complement are shown.  Parallel lines drawn in black denote a perfect matching and in white denote the absence of one.  Shaded areas indicate that all possible edges are present.  \label{fig:R4m}}
\end{figure}
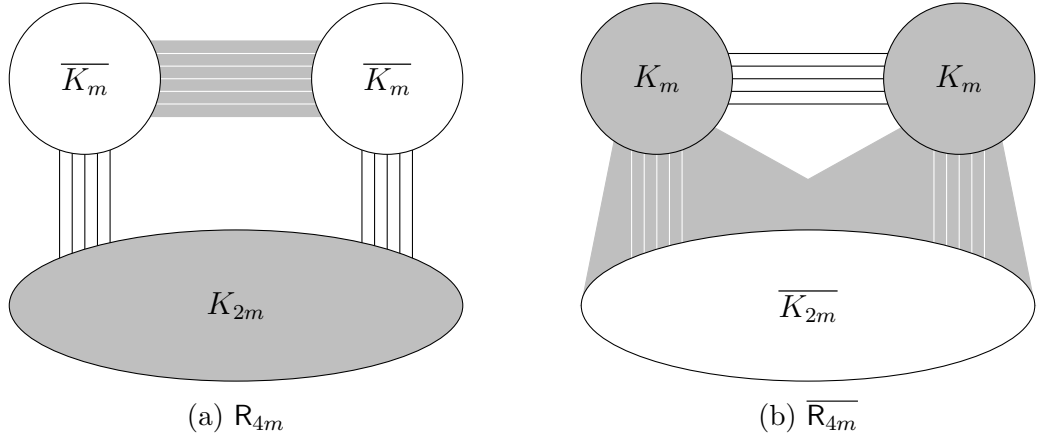

 \begin{ex}\label{Rn}
     We define a graph $G=\rn_{4m}$ on $4m$ vertices as  shown  in Figure \ref{fig:R4m}(a).   Begin with three sets $X$, $Y$, and $W$ such that $|X|= |Y|=m$, and $|W|=2m$ and let $V(G)=X\cup Y\cup W$. Then add edges to construct $G$ as follows: Make $G[W] \cong K_{2m}$ and add a perfect matching between $X \cup Y$ and $W$. Finally, include all edges between $X$ and $Y$ except for a perfect matching. Now $G$ and $\ol G$ are both connected. Note that $W$ is a  standard forcing set for $G$ and $X \cup W$ is a  standard forcing set for $\ol G$, each with 
     propagation time one. Hence,  $\thza(G) \le  2m $ and $\thza(\ol G) \le  3m $. 
 \end{ex}

We use the same partition strategy for  the next two theorems and much of the proofs are identical, so we first define the partition and prove  some basic properties in a lemma.  

Let $G_1$ and $G_2$ be disjoint graphs of equal order and let $M$ be a matching between $V(G_1)$ and $V(G_2)$ that saturates all vertices. Then the  \emph{$M$-sum of $G_1$ and $G_2$}, denoted by $G_1M^+G_2$, is the graph with $V(G_1M^+G_2)=V(G_1)\cup V(G_2)$ and $E(G_1M^+G_2)=E(G_1)\cup E(G_2)\cup M$. A graph of the form $G_1M^+G_2$ is also called a \emph{matched-sum graph}.  The graph $\rn_{4m}$ in Example \ref{Rn} is an example of a matched-sum graph.

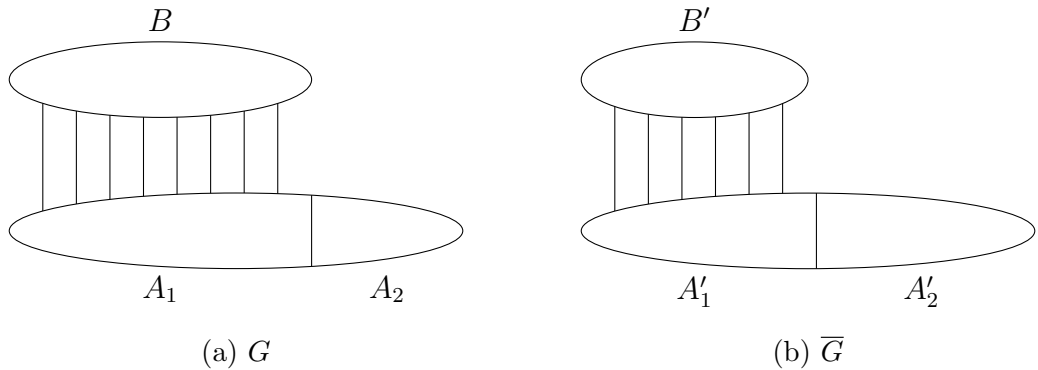
\begin{figure}[!h]
    \centering
    \begin{subfigure}[b]{0.45\textwidth}
    \centering
    \raisebox{1.8pt}{\begin{tikzpicture}
    \foreach \x in {1,...,8}{
        \draw (-3+0.444*\x,2) -- (-3+0.444*\x,0);
    }
    \draw[fill=white] (0,0) ellipse (3 and 0.5);
    \draw (1,{sqrt(2)/3}) -- (1,{-sqrt(2)/3});
    \draw[fill=white] (-1,2) ellipse (2 and 0.5);
    \draw (-1,2.5) node[anchor=south]{$B$};
    \draw (-1,-0.45) node[anchor=north]{$A_1$};
    \draw (2,-0.45) node[anchor=north]{$A_2$};
    \end{tikzpicture}}
    \caption{$G$}
    \end{subfigure}
    \begin{subfigure}[b]{0.45\textwidth}
    \centering
    \begin{tikzpicture}
    \foreach \x in {1,...,6}{
        \draw (-3+0.444*\x,2) -- (-3+0.444*\x,0);
    }
    \draw[fill=white] (0,0) ellipse (3 and 0.5);
    \draw (1/9,0.5) -- (1/9,-0.5);
    \draw[fill=white] (-1.5,2) ellipse (1.5 and 0.5);
    \draw (-1.5,2.5) node[anchor=south]{$B'$};
    \draw (-1.5,-0.45) node[anchor=north]{$A_1'$};
    \draw (1.5,-0.45) node[anchor=north]{$A_2'$};
    \end{tikzpicture}
    \caption{$\ol G$}
    \end{subfigure}
    \caption{ $M^+$-partitions of the graphs  $G$ and $\ol G$ depicted as in the proofs of Lemma \ref{NGlower-Zast-lem} and   Theorems \ref{NG-Z-con-lower}  and \ref{NG-Z-prod-lower}. In (a) there may be edges in  $G[A_1]$, $G[A_2]$, and $G[B]$, and between $A_2$ and $A_1\cup B$; similar comments apply to (b).}
    \label{fig:R-lower}
\end{figure}

Let $H$ be a graph. A partition of $V(H)$ into a triple of sets $(A_1,A_2,B)$  is an \emph{$M^+$-partition} if  $H[A_1\cup B]\cong H[A_1]M^+H[B]$ 
(see Figure \ref{fig:R-lower}). 

\begin{lem}\label{NGlower-Zast-lem} Let $G$ be a graph of order $n$.   There exist $A, A' \subseteq V(G)$ such that $\thza(G)=\thza(G;A)$ and $\thza(\ol G)=\thza(\ol G; A')$,  and there is a  forcing process for each of these sets that forces in one round.  Define $A_1\subseteq A$ and $A'_1\subseteq A'$ to be the  sets of vertices that perform forces  in these processes.  
 Let $A_2=A\setminus A_1$ and $A'_2=A'\setminus A'_1$, and let $B = V(G)\setminus A$ and $B' = V(\ol G) \setminus A'$.   Then: 
\ben[$(1)$]
\item\label{NGlower-Zast-lem:c1} $(A_1,A_2,B)$ is an $M^+$-partition and $(A'_1,A'_2,B')$ is an $M^+$-partition of $\Gc$.
\item \label{Eq:PartitionSum} 
$|A_1|+|A_2|+|B|=n \hspace{.3cm} \text{ and } \hspace{.3cm} |A'_1|+|A'_2|+|B'|=n.$
\item\label{Eq:A1andB} $|A_1|=|B| \hspace{.3cm}\text{ and } \hspace{.3cm} |A'_1|=|B'|.$
\item \label{eq:thza+thzacpart}
   $\thza(G)=n-|A_1| = \frac{n+|A_2|}{2} \hspace{.3cm} \mbox{ and }\hspace{.3cm} \thza(\Gc)= n - |A'_1| = \frac{n+|A'_2|}{2} . $
   \item\label{NGlower-Zast-lem:c5} The edges between $A_1$ and $B$ in $\ol G[A_1\cup B]$ form a complete bipartite graph with a perfect matching removed, and  the edges between $A_1'$ and $B'$ in $G[ A'_1\cup B']$ form another complete bipartite graph with a perfect matching removed. 
   \item\label{switching} 
     $B\cup A_2$ is a forcing set for $G$ with $\thza(G)=\thza(G; B\cup A_2)$, and $B$ is the set of filled vertices that force $A_1$ in a single round. When we consider the forcing set $B\cup A_2$ instead of $A_1 \cup A_2$, we refer to this as \emph{switching} between $A_1$ and $B$. 
\item\label{Eq:A'capA}\label{Eq:A'capB}
  $   A'_1 \cap A_1\neq \emptyset $ implies $ |B\cap B'|\leq 2 
$, and $
    A'_1 \cap B\neq \emptyset $ implies $ |A_1\cap B'|\leq 2.
$
\item\label{onempty}  $(A'_1 \cap A_1= \emptyset$ and $A'_1\cap B \neq \emptyset)$ implies  $|A_1|-2\leq |A'_2|$.
\item \label{bothnonempty}   $(A'_1 \cap A_1\neq \emptyset$ and $A'_1 \cap B\neq \emptyset )$ implies $|A_1|+|B|-4\leq \thza(\ol G)$. 

 \een

\end{lem}

\bpf 
   Sets $A$ and $A'$  that realize the throttling numbers and  are able to  force all vertices of $G$ and $\Gc$ in a single round exist by  Theorem \ref{t:NG-NCproductZset}.   This implies \eqref{NGlower-Zast-lem:c1} and \eqref{Eq:PartitionSum}.  The  perfect matchings required by a $M^+$-partition establish \eqref{Eq:A1andB},  \eqref{NGlower-Zast-lem:c5}, and that
  $G[A_1\cup B]=G[A_1]M^+G[B]$ and $\ol G[ A'_1\cup B']=\ol G[A'_1]M^+ \ol G[B']$.  Parts \eqref{Eq:PartitionSum} and \eqref{Eq:A1andB} imply \eqref{eq:thza+thzacpart}.  Part \eqref{switching} follows by  reversing the forcing process.

 For Part \eqref{Eq:A'capA}, assume $A'_1 \cap A_1$ is nonempty with a vertex $x\in A'_1 \cap A_1$. In $G$, $x$ forces some vertex $y\in B$; in $\ol G$, $x$ forces  some vertex $y'\in B'$. Suppose there exists a vertex $z\in(B\cap B')\setminus\{y,y'\}$. By the previous remarks, since $xy$ is an edge in $G$, $xz$ must be an edge in $\ol G$. Similarly, since $xy'$ is an edge in $\ol G$, $xz$ must be an edge in $G$, a contradiction. Therefore, $(B\cap B')\setminus\{y,y'\}=\emptyset$, establishing the first inequality in \eqref{Eq:A'capA}. By  switching between $A_1$ and $B$, the second inequality is established.

 For Part \eqref{onempty}, assume the hypotheses.   Since $
    A'_1 \cap B\neq \emptyset $,  \eqref{Eq:A'capB} implies $|A_1\cap B'|\leq 2$. Since $A'_1 \cap A_1= \emptyset$, 
    we have  $A_1\setminus (A_1\cap B')=A_1\setminus B'\subseteq  A'_2$. Thus  $|A_1|-2\le |A_1|-|A_1\cap B'|\le   |A'_2|$.

 For Part \eqref{bothnonempty},  suppose $A'_1 \cap A_1\neq \emptyset$ and $A'_1 \cap B\neq \emptyset$. Since $(A_1\cup B)\setminus ((B\cap B')\cup(A_1\cap B'))=(A_1\cup B)\setminus B' \subseteq \ol{B'}=A'_1\cup A'_2$, \eqref{Eq:A'capA}  gives us $|A_1|+|B|-4\leq |A'_1|+|A'_2|=\thza(\ol G)$. 
\epf

\begin{thm}\label{NG-Z-con-lower}
Let $G$ be a graph of order $n\ge 2$.  Then $\frac{5n}{4} - 2 \le \thza(G)+ \thza(\ol G)$ and this bound 
is almost sharp even when both $G$ and $\ol G$ are connected.
\end{thm}

\bpf
We begin by partitioning $V(G)$ into sets $A_1,A_2$ and $B$ and partitioning $V(\ol G)$ into sets $A'_1$, $ A'_2$ and $B'$ as in Lemma \ref{NGlower-Zast-lem}.
For each of the following cases, we assume without loss of generality that $\thza(G)\leq \thza(\ol G)$.

Case 1: Suppose $A'_1 \cap A_1\neq \emptyset$ and $A'_1 \cap B\neq \emptyset$. 
 By Lemma \ref{NGlower-Zast-lem}\eqref{bothnonempty}, $|A_1|+|B|-4\leq \thza(\ol G)$. If $\thza(G)\geq \frac{5n}{8}-1$,  then $\frac{5n}{4} - 2 \le \thza(G)+ \thza(\ol G)$. Hence, we can assume that $\thza(G)< \frac{5n}{8}-1$. By Lemma \ref{NGlower-Zast-lem}\eqref{eq:thza+thzacpart}, 
$n-|A_1|<\frac{5n}{8}-1$ and so $\frac{3n}{8}+1<|A_1|$.  
Additionally, we know by Theorem \ref{t:NG-NCproductZset} that  $\thza(G)\geq \frac{n}{2}$. Therefore, $\thza(G)+ \thza(\ol G)\geq \frac{n}{2} + |A_1|+|B|-4=\frac{n}{2}+2|A_1|-4>\frac{n}{2}+2(\frac{3n}{8}+1)-4=\frac{5n}{4} - 2$.

Case 2: Suppose $A'_1 \cap A_1= \emptyset$ and $A'_1 \cap B\neq \emptyset$.  By Lemma \ref{NGlower-Zast-lem}\eqref{onempty}, $|A_1|-2\leq |A'_2|$.  
By Lemma \ref{NGlower-Zast-lem}\eqref{eq:thza+thzacpart}, 
$\thza(G)+ \thza(\ol G)= n-|A_1|+\frac{n+|A'_2|}{2} =\frac{3n}{2}-|A_1|+\frac{|A'_2|}{2}\geq \frac{3n}{2}-|A_1|+\frac{|A_1|-2}{2}=\frac{3n}{2}-\frac{|A_1|}{2}-1\geq \frac{3n}{2}-\frac{n}{4}-1=\frac{5n}{4}-1$.

Case 3: Suppose $A'_1 \cap A_1\neq \emptyset$ and $A'_1 \cap B= \emptyset$. By switching between $A_1$ and $B$, this case is equivalent to Case 2.

Case 4: Suppose $A'_1 \cap A_1= \emptyset$ and $A'_1 \cap B= \emptyset$.  If $\thza(G) \geq \frac{5n}{8}$, the result follows. So we may assume that $\frac{5n}{8} > \thza(G)= n-|A_1|$. 
This implies $\frac{3n}{8} < |A_1|$, and hence $|A_2| < \frac{5n}{8} - |A_1| < \frac{n}{4}$.  Consequently $A'_1 \subseteq \ol{A_1 \cup B} = A_2$ implies $|A'_1| < \frac{n}{4}$.  By Lemma \ref{NGlower-Zast-lem}\eqref{Eq:PartitionSum}  and \eqref{Eq:A1andB}, $\thza(G) + \thza(\ol G) = |A_1|+|A_2|+|A'_1|+|A'_2| \geq \frac{n}{2} + n - |A'_1| > \frac{5n}{4} $.

In each of Cases 1, 2, 3, and 4 we arrive at the conclusion that $\frac{5n}{4} - 2 \le \thza(G)+ \thza(\ol G)$. 
 Example \ref{Rn} shows this bound is almost sharp because $\thza(\rn_{4m}) + \thza(\ol {\rn_{4m}}) \leq 2m + 3m = 5|V(\rn_{4m})|/4$.
\epf

\begin{thm}\label{NG-Z-prod-lower}
Let $G$ be a graph of order $n\ge 4$.  Then $\frac{3n^2}{8} - O(n) \le \thza(G) \cdot \thza(\ol G)$ and this bound is tight even when both $G$ and $\ol G$ are connected.
\end{thm}

\bpf 
We begin by partitioning $V(G)$ into sets $A_1,A_2$ and $B$ and partitioning $V(\ol G)$ into sets $A'_1$, $ A'_2$ and $B'$ as in Lemma \ref{NGlower-Zast-lem}. For each of the following cases, we assume that we have labeled $G$ and $\ol G$ such that $\thza(G)\leq \thza(\ol G)$.

Case 1: Suppose $A'_1 \cap A_1\neq \emptyset$ and $A'_1 \cap B\neq \emptyset$. 
 By Lemma \ref{NGlower-Zast-lem}\eqref{bothnonempty}, $|A_1|+|B|-4\leq \thza(\ol G)$. 
If $\thza(G)\geq \frac{\sqrt{6}n}{4}$,  then $\frac{3n^2}{8} \le \thza(G)\cdot\thza(\ol G)$. Hence, we can assume that $\thza(G)< \frac{\sqrt{6}n}{4}$. By   Lemma \ref{NGlower-Zast-lem}\eqref{eq:thza+thzacpart} $n-\frac{\sqrt{6}n}{4}<|A_1|$.  Again by Theorem \ref{t:NG-NCproductZset}, we know $\thza(G)\geq \frac{n}{2}$. 
We conclude that \[\thza(G)\cdot\thza(\ol G)\geq \lp\frac{n}{2}\rp(2|A_1|-4)>\lp\frac{n}{2}\rp\lp 2\lp n-\frac{\sqrt{6}n}{4}\rp-4\rp=\frac{4-\sqrt{6}}{4}n^2-2n > \frac 3 8 n^2-2n.\]

Case 2: Suppose $A'_1 \cap A_1= \emptyset$ and $A'_1 \cap B\neq \emptyset$. 
 By Lemma \ref{NGlower-Zast-lem}\eqref{onempty}, $|A_1|-2\leq |A'_2|$. 
By  Lemma \ref{NGlower-Zast-lem}\eqref{eq:thza+thzacpart}, \[\thza(G)\cdot\thza(\ol G)=(n-|A_1|)\lp\frac{n+|A'_2|}{2}\rp
\geq (n-|A_1|)\lp\frac{n}{2}+\frac{|A_1|}{2}-1\rp=\frac{n^2}{2}-n
-\frac{|A_1|^2}{2}+|A_1|.\]

\noindent Since $|A_1|\leq \frac{n}{2}$, for $n\geq 4$ we have $-\frac{|A_1|^2}{2}+|A_1|\geq -\frac{n^2}{8}+\frac{n}{2}$. Hence, $\thza(G)\cdot\thza(\ol G)\geq \frac{3n^2}{8}-\frac{n}{2}$.

Case 3: Suppose $A'_1 \cap A_1\neq \emptyset$ and $A'_1 \cap B= \emptyset$. By switching between $A_1$ and $B$, this case is equivalent to Case 2.

Case 4: Suppose $A'_1 \cap A_1= \emptyset$ and $A'_1 \cap B= \emptyset$. Then $A'_1 \subseteq \ol{A_1 \cup B} = A_2$. If $\thza(G) \geq \frac{\sqrt{6}n}{4}$, the result follows. So we may assume that $\frac{\sqrt{6}n}{4} > \thza(G)=|A_1|+|A_2|$.  As in Case 1,  
Lemma \ref{NGlower-Zast-lem}\eqref{eq:thza+thzacpart} implies
$n-\frac{\sqrt{6}n}{4} < |A_1|$, and hence $|A_2| = n-2|A_1|< \frac{(\sqrt{6}-2)n}{2}$. Then $A'_1 \subseteq A_2$ implies $|A'_1|<\frac{(\sqrt{6}-2)n}{2}$. 
As before, 
\[\thza(G)\cdot\thza(\ol G) \geq 
\lp\frac{n}{2}\rp\lp n - |A'_1|\rp
> \lp\frac{n}{2}\rp\lp n - \frac{(\sqrt{6}-2)n}{2}\rp
=\frac{4-\sqrt{6}}{4}n^2.\] 

In each of Cases 1, 2, 3, and 4 we arrive at the conclusion that $\frac{3n^2}{8} - O(n)  \le \thza(G)\cdot \thza(\ol G)$. 
 Example \ref{Rn} shows this bound is tight because 
\[\thza(\rn_{4m})\cdot\thza(\ol{\rn_{4m}}) \leq (2m)(3m) = \frac{3(4m)^2}{8} = \frac{3|V(\rn_{4m})|^2}{8}.\qedhere\]
\epf


\section{PSD zero forcing}  \label{s:PSD-NG-th}
In this section we establish Nordhaus-Gaddum sum and product bounds for $\thpa$ and $\thpx$.   We improve the known NG sum lower bound for $\thp$ and establish NG product bounds for this parameter.
We also provide NG sum and product bounds for all three parameters when  both $G$ and $\Gc$ must be connected.  The next result describes previously known Nordhaus-Gaddum bounds for sum throttling of PSD forcing.  

\begin{thm}\label{NG-Zp}{\rm \cite[Theorem 11.39]{HLSbook}} Let $G$ be a graph  of order $n$.  
 Then  $n+o(n)\le \thp(G)+\thp(\Gc)\le 2n$. If each of $G$ and $\Gc$ is required to have an edge, then $n+o(n)\le \thp(G)+\thp(\Gc)\le 2n-1$.  The upper  bounds are sharp and the lower bounds are tight.
\end{thm}

 We begin with a lemma on NG product lower bounds that will be used in both the general case and the case in which both graphs must be connected.    Recall the straightforward fact that
 $\Zp(G)\ge \frac n 2$ implies $\thpx(G)=n$ for any graph $G$ of order $n$.

  \begin{lem}  \label{Lem:NG-Zp-prodth-lower-arb}
Let $G$ be a graph of order $n\ge 6 $ such that $\Zp(G)\le\Zp(\Gc)$. 
Then:
 \ben[$(1)$] 
 \item\label{NGZpthbdsL}  $2n \le \thp(G)\cdot \thp(\ol G) $,  and $3n-9 \le \thp(G)\cdot \thp(\ol G) $ if  $G\not\cong K_{1,n-1}$. 
\item\label{NGZpprodthabdsL} $  n-1\le \thpa(G)\cdot \thpa(\ol G) $, and $ 2n-8\le \thpa(G)\cdot \thpa(\ol G) $ if $G\not\cong K_{1,n-1}$. 
\item\label{NGZpprodthxbdsL}  $  2n\le \thpx(G)\cdot \thpx(\ol G)  $, and    $  3n\le \thpx(G)\cdot \thpx(\ol G)  $ if   and $G\not\cong K_{1,n-1}$. 
\een
Each of the first lower bounds is sharp, realized by $K_{1,n-1}$.  \end{lem}

 \bpf  
We examine three cases based on $\Zp(G)$ and $\ptp(G)$.
 Consider first the case in which $\Zp(G)=1$ and $\ptp(G)=1$.  This is equivalent to  $G\cong K_{1,n-1}$; note that  $\ol{K_{1,n-1}}\cong K_{n-1}\du K_1$, so $\Zp(\ol{K_{1,n-1}})=n-1$ and $\ptp(\ol{K_{1,n-1}})=1$.  Thus
\[\thp(K_{1,n-1})=2, \qquad\thpa(K_{1,n-1})=1, \qquad\thpx(K_{1,n-1})=2.\]
\[\thp(\ol{K_{1,n-1}})=n, \qquad\thpa(\ol{K_{1,n-1}})=n-1, \qquad\thpx(\ol{K_{1,n-1}})=n, \]
which shows that the  first lower bound in each statement (once established)  is sharp.  
 For the remaining two cases, we will use the fact that    $\Zp(G)+\Zp(\Gc)\ge n-2$ \cite[Theorem 9.60]{HLSbook}.
 
  For the second case,  suppose that $\Zp(G)=1$ and $\ptp(G)\ge 2$.  Then $\thpx(G)\ge 3$ and $\Zp(\Gc)\ge n-3$. 
  Since $n\ge 6$, $\thpx(\Gc)=n$, so $\thpx(G)\cdot\thpx(\Gc)\ge 3n$. 
Also $\thpa(G)\ge 2$ and   $\thpa(\Gc)\ge n-3$, so $\thpa(G)\cdot\thpa(\Gc)\ge 2n-6\ge n-1$  since $n\ge 6$ (and $2n-6\ge 2n-8$ for all $n$). 
 Furthermore, $\thp(G)\ge 3$,  and $\thp(\Gc)\ge n-2$. Thus $\thp(G)\cdot\thp(\Gc)\ge 3n-6\ge 2n$  since $n\ge 6$ (and $3n-6\ge 3n-9$ for all $n$).  

 For the third case, assume $\Zp(G)  \ge 2$.       
 First we consider $\thpx$ and note that $ \thpx(G)\ge 4$. It has  been verified in \cite{sage} that $\thpx(G)\cdot \thpx(\Gc)\ge 3n$ for $ n=6$ (since $G\not\cong K_{1,n-1}$).  Since  $n=7$ implies  $\Zp(\Gc)\ge 3$ and thus $\thpx(G)\cdot\thpx(\Gc)\ge 24>21$, assume $n\ge 8$.  Then $\thpx(G)\cdot\thpx(\Gc) \geq 4\lp\frac{n}{2} - 1\rp(2) = 4n - 8  \ge 3n$, so  both lower bounds are established and this  completes  the proof for $\thpx$.

Next we consider  $\thpa$ with $\Zp(G) = z \ge 2$. 
First assume $2\le z\le \frac n 2 -1$. To show that  $\thpa(G)\cdot \thpa(\Gc)\ge 2n-8$  for $2\le z\le \frac n 2-1 $, it suffices to show that $0\le g(z)= z(n-2-z)-(2n-8) =-z^2+(n-2)z-2n+8$.  This is verified by evaluating $g(2)=0$ and $g\lp\frac n 2 -1\rp=\frac{1}{4}(n-6)^2\ge 0$.  Since the last inequality also implies $\lp \frac n 2-1\rp^2\ge 2n-8$, we have   $\thpa(G)\cdot \thpa(\Gc)\ge 2n-8$ for $z\ge \frac n 2-1$  because $\Zp(\Gc)\ge z$.   It has  been verified in \cite{sage} that $\thpa(G)\cdot \thpa(\Gc)\ge n-1$ for $ n=6$. 
Since $n-1\le 2n-8$ for $n\ge 7$, this completes the proof for $\thpa$.

Finally, we consider  $\thp$ with $\Zp(G) = z \ge 2$. First assume $2\le z\le \frac n 2 -1$.  To show $\thp(G)\cdot\thp(\Gc)\ge 3n-9$, it suffices to show that $0\le f(z)= (z+1)(n-2-z+1)-(3n-9) =-z^2+(n-2)z-2n+8$ for $2\le z\le \frac n 2-1 $.  This is verified by evaluating $f(2)=0$ and $f\lp\frac n 2 -1\rp=\frac{1}{4}(n-6)^2\ge 0$.  For  $z\ge \frac n 2-1$,   $\thp(G)\cdot \thp(\Gc)\ge (z+1)^2 \ge \lp \frac n 2\rp^2\ge 3n-9$.   It has also been verified in \cite{sage} that $\thp(G)\cdot \thp(\Gc)\ge 2n$ for $n=6,7,8$. Since $2n \le 3n-9$ for $n\ge 9$, this completes the proof for $\thp$.  
 \epf

\begin{rem}  \label{rem:NG-Zp-prodth-lower-arb}
Lemma \ref{Lem:NG-Zp-prodth-lower-arb} does not address graphs of orders $n=4,5$ because there are some some special cases for these orders.  Here we list sharp bounds when $n=4$ and $5$, all determined in \cite{sage}.  
 \ben[$(1)$] 
 \item  $8 \le \thp(G)\cdot \thp(\ol G) $  for $n=4$ and $9 \le \thp(G)\cdot \thp(\ol G) $  for $n=  5$.
\item $  3\le \thpa(G)\cdot \thpa(\ol G) $ for $n=4$ and  $   4\le \thpa(G)\cdot \thpa(\ol G) $ for $n=5$. 
\item  $  8\le \thpx(G)\cdot \thpx(\ol G)  $ for $n=4$ and $  10\le \thpx(G)\cdot \thpx(\ol G)  $ for $n=5$.
\een
Each of the general lower bounds is sharp,  realized by $G\cong K_{1,n-1}$  except for $\thp$ and $n= 5$. 
  The graphs of order five having $\thp(G)=\thp(\Gc)=3$ are $ C_5,$ the Bull Graph, and $P_5$ and its complement the House Graph (pictured in \cite{sage}). 
  \end{rem}

\begin{cor}
 \label{NG-Zp-prodth-both-arb}
Let $G$ be a graph of order $n\ge 4$. 
Then:
 \ben[$(1)$] 
\item\label{NGZpthbds} 
$  n\le \thp(G)+ \thp(\ol G) \le 2n$;\\  for $n\ne  5$, $2n \le \thp(G)\cdot \thp(\ol G) \le n^2$   (for $n=  5$, $9 \le \thp(G)\cdot \thp(\ol G) \le  25$).

\item\label{NGZpprodthabds} $  n-2\le \thpa(G)+ \thpa(\ol G) \le 2n-2$ and $  n-1\le \thpa(G)\cdot \thpa(\ol G) \le (n-1)^2$. 
\item\label{NGZpprodthxbds} $  n\le \thpx(G)+ \thpx(\ol G)  \le 2n$ and $  2n\le \thpx(G)\cdot \thpx(\ol G)  \le n^2$.  
\een
All NG upper bounds are sharp; the NG product lower bounds are sharp; the NG sum lower bounds are almost sharp.
\end{cor}
 \bpf  The NG upper bounds are all immediate by Remark \ref{NG-u}. The graph $K_n$ shows these bounds are sharp for $\thp$ and $\thpx$; for $\thpa$ use $K_n-e$ (to guarantee both graphs have an edge). 

The NG product lower bounds  and their sharpness were established in Lemma \ref{Lem:NG-Zp-prodth-lower-arb}   (and Remark \ref{rem:NG-Zp-prodth-lower-arb} for $n=4,5$). The  NG sum lower bounds follow from $\Zp(G)+\Zp(\Gc)\ge n-2$ \cite[Theorem 9.60]{HLSbook}  and Remark \ref{NG-u}; $K_{1,n-1}$ shows these bounds are almost sharp (see values of throttling numbers in the proof of Lemma \ref{Lem:NG-Zp-prodth-lower-arb}). 
\epf

Next we turn our attention to the case when both graphs must be connected.  The example provided in \cite[Proposition 11.39]{HLSbook} for the NG sum upper bound for $\thp$   
is $K_n$ (or if each graph must have an edge then  $K_n-e$), which cannot be used when both $G$ and $\ol G$ must be connected. The following example is useful for  the NG sum and product upper bounds for $\thp$, $\thpa$,   and $\thpx$ when both graphs must be connected.  

\begin{ex}\label{ex:Fc} For $c\ge 3$, construct a graph $\fc_c$ of order $c^2$  by adding a cycle of $ c$ edges to the disjoint union of $c$ copies of $K_c$ (numbered $0,\dots, c-1)$. Label the set of  vertices of the $i$th copy of $K_c$ by $V_i=\{x_i,u_{i,1},\dots,u_{i,c-1}\}$ where $x_i$ is the cycle vertex, and denote by $E_i$ the edges  within  the $i$th copy of $K_c$.  Then $V(\fc_c)=\bigcup_{i=0}^{c-1} V_i$ and $E(\fc_c)=\bigcup_{i=0}^{c-1} E_i\cup \{x_0x_1,x_1x_2,  \dots, x_{c-1}x_0\}$.  The graph $\fc_6$ is shown in Figure \ref{fig:F6}.  
Let $X=\{x_0,\dots,x_{c-1}\}$ and $U_i=\{u_{i,1},\dots,u_{i,c-1}\}$ for $i=0,\dots,c-1$.
Observe that $\ol{\fc_c}$ can be constructed from $K\!\!\underbrace{_{c,c,\dots,c}}_{c\ {\rm times}}$ by removing a $c$-cycle of edges.  
\end{ex} 

\begin{figure}[!h]
    \centering
    \begin{tikzpicture}
        \foreach \t in {0,...,5}{
            \draw (30+60*\t:2) -- (90+60*\t:2);
            \foreach \a in {0,...,5}{
                \foreach \b in {\a,...,4}{
                    \draw ($(30+60*\t:2.75) + (30+60*\a:0.75)$) -- ($(30+60*\t:2.75) + (90+60*\b:0.75)$);
                }
            }
        }
        \foreach \t in {0,...,5}{
            \draw (30+60*\t:1.65) node{$x_{\t }$};
            \foreach \a in {0,...,5}{
                \draw[fill=white] (30+60*\t:2.75) + (30+60*\a:0.75) circle (0.1);
            }
        }
    \end{tikzpicture}
    \caption{The graph $\fc_6$. }
    \label{fig:F6}
\end{figure}

In the next proof we  use PSD forts, which were  introduced by Smith, Mikesell, and Hicks in \cite{SMH19}. Let $F\subseteq V(G)$ be nonempty and partition $F=F_1\du\dots\du F_k$ where $G[F_i]$ are the connected components of $G[F]$.  Then $F$  is a \emph{PSD fort} of $G$ if and only if for every $v\in V(G)\setminus F$, $|N_G(v)\cap F_i|\ne 1$ 
for $i=1,\dots, k$.  Furthermore, they established the following result.

\begin{thm} \label{t:PSDfort}  {\rm \cite{SMH19}}
For a graph $G$, $S\subseteq V(G)$ is a PSD forcing set if and only if $S\cap F\ne \emptyset$ for every PSD fort $F$ of $G$.
\end{thm}

\begin{thm}\label{p:Fc} Let $\fc_c$ be the graph defined in Example \ref{ex:Fc}.  Then 
\[ \Zp(\fc_c)=c^2-2c+2, \qquad  \Zp(\ol{\fc_c})=c^2-c,\]
\[ \thp(\fc_c)\ge c^2-2c+3, \qquad \thp(\ol{\fc_c})=c^2-c+1,\]
\[ \thpa(\fc_c)=\thpa(\ol{\fc_c})=c^2-c, \qquad \text{and} \qquad \thpx(\fc_c)=\thpx(\ol{\fc_c})=c^2.\]
\end{thm}

\bpf 
First we exhibit minimum PSD forcing sets.
The set $B_0=V(\fc_c)\setminus (X\cup\{u_{i,1}:i=0,\dots,c-1\})\cup \{x_i,x_k\}$ with $i\ne k$ is a PSD forcing set for $\fc_c$:  Indeed, the deletion of $x_i, x_k$ disconnects the graph into three or four connected components, so that $x_i$ can PSD force $u_{i,1}, x_{i+1}, x_{i-1}$ (if $i = k\pm 1$, then the force of an already filled vertex does not take place), where all arithmetic in the indices is considered modulo $c$; the initial forces for $x_k$ occur similarly, and the process can then be repeated. Thus $\Zp(\fc_c)\le c^2-2c+2$. 
The set $B'_0=V(\ol{\fc_c})\setminus V_1$ is a PSD forcing set for $\ol{\fc_c}$ because $ \ol{\fc_c}[V(\ol{\fc_c})\setminus B_0'] = \ol{\fc_c}[V_1]$ is a set of isolated vertices. Thus $\Zp(\ol{\fc_c})\le c^2-c$.

To see that $\Zp(\fc_c)\ge c^2-2c+2$, let $B\subseteq V(\fc_c)$ be a PSD forcing set. Then $|U_i\cap B|\ge c-2$ for every $i$, thus $|B|\geq c(c-2)=c^2-2c$. If $|B|\leq c^2-2c+1$, then $|X\cap B|\leq1$ and all vertices in $X\setminus B$ can never by PSD forced. Hence, $|B|\geq c^2-2c+2$.

To see that $\Zp(\ol{\fc_c})\ge c^2-c$, we first note that this has been verified for  $c=3,4,5$ \cite{sage}, so we assume $c\ge 6$. 
Let $B\subseteq V(\ol{\fc_c})$ be a PSD forcing set, let $W=V(\ol{\fc_c})\setminus B$, and suppose $|W|\ge c+1$.  We consider several cases, based on the size of the largest intersection of $W$ and some $V_i$, and show that in each case $B$ is not a PSD forcing set, producing a contradiction. 
Recall that for any $p\ne q$,   at most one edge between vertices in $V_p$ and $V_q$ is not in $\ol{\fc_c}$.

 Suppose first that $|W\cap V_p|=c$ for some $p$.  Either $\ol{\fc_c}[W]$ is connected or $W$ can be partitioned as $(W\setminus\{x_p\}) \du \{x_p\}$ to yield connected components.  Note that every vertex of $B$ is adjacent to all $c-1\ge 2$ vertices of $U_p\subseteq W\setminus\{x_p\}$.  In the first case, this immediately implies that $W$ is a fort of $\ol{\fc_c}$ and $B$ is not a PSD forcing set.  In the second case, since $x_p$ has no neighbors in $W\setminus \{x_p\}$, we see that $W\setminus\{x_p\}$ is a fort and $B$ is not a PSD forcing set.

Suppose next that $|W\cap V_p|=c-1$ for some $p$. Then there are at least two unfilled vertices $v,v'\in V(G)\setminus V_p$. Furthermore, $|W\cap U_p|\geq2$. As long as at least one of $v$ and $v'$ remains unfilled, all the unfilled vertices of $U_p$ are in one connected component of the subgraph induced by the currently unfilled vertices, so no vertex can force any unfilled vertex in $U_p$. The only vertex that can possibly force $v$ or $v'$ is $x_p$ (in the case that $|\{v,v'\}\cap\{x_{p+1},x_{p-1}\}|=1$), but at most one of $v$ and $v'$ can be forced and no further PSD forces can occur. Thus $B$ is not a PSD forcing set.

Subsequently, suppose that $3\leq|W\cap V_p|\leq c-2$ for some $p$. Then there are at least three unfilled vertices $v,v',v''\in V(G)\setminus V_p$. With a similar observation as the previous case, all of $v$, $v'$, and $v''$ need to be filled before any of the unfilled vertices in $U_p$ can be filled. The only vertex that can possibly force $v$, $v'$, or $v''$ is $x_p$ (in the case that $|\{v,v',v''\}\cap\{x_{p+1},x_{p-1}\}|=2$), but at most one of $v$, $v'$, and $v''$ can be forced and no further PSD forces can occur. Thus $B$ is not a PSD forcing set.

Finally suppose that $|W\cap V_i|\leq 2$ for every $i$.  Since $|W|\geq c+1$, there exists some $p$ such that $|W\cap V_p| = 2$; let $W\cap V_p = \{v_p, u_p\}$ where $u_p \in U_p$.  Then $u_p$ is adjacent to every other vertex in $W$ except $v_p$.  Since $c+1\geq 7$, $|\{i: W\cap V_i \neq \emptyset\}| \geq 4$, so $v_p$ has at least one neighbor in $W$.  Therefore $\ol{\fc_c}[W]$ is connected.  Given any $v\in B$, suppose $v\in V_j$ for some $j$; since $|W\cap V_j|\leq 2$ and $v$ has at most two non-neighbors in $V(\ol{\fc_c})\setminus V_j$, we have $|N(v)\cap W| \geq c+1-2-2 \geq 3$.  So $W$ is a PSD fort, and $B$ is not a PSD forcing set.

We complete the proof by establishing the statements about throttling numbers. Note that $\thpx(\fc_c)= \thpx(\ol{\fc_c})=c^2$ is immediate from  $\Zp(\fc_c)= c^2-2c+2$  and $\Zp(\ol{\fc_c})= c^2-c$.  Furthermore, $\thp(\fc_c)\ge c^2-2c+3$ is immediate from  $\Zp(\fc_c)= c^2-2c+2$  and Remark \ref{NG-u2}.  Since $\ptp(\ol{\fc_c},c^2-c)=1$,  $\thp(\ol{\fc_c})=c^2-c+1$ and $\thpa(\ol{\fc_c})=c^2-c$. 

Lastly we show that $\thpa(\fc_c) = c^2-c$.  Since $\ptp(\fc_c, c^2-c) = 1$, $\thpa(\fc_c) \leq c^2-c$ is immediate.  Suppose $B$ is a set such that $\thpa(\fc_c; B) = \thpa(\fc_c) < c^2-c$.  Since $|B|\geq c^2 - 2c + 2$ and $t(c^2 - 2c + 2) > c^2 - c$ for all $t\geq 2$, we must have $\ptp(\fc_c; B) = 1$ and  $|B| < c^2 - c$. 
So there are at least $c+1$ unfilled vertices, which implies $|B\cap V_i| \leq c-2$ for some $i$.  Then there is a vertex in $U_i$ which cannot be forced in the first round, and the contradiction implies $\thpa(\fc_c) \geq c^2-c$.  
\epf

We can now establish NG sum and product bounds, with specified sharpness or tightness, for all three types of throttling of PSD forcing when both $G$ and $\Gc$ are connected.
Note that  Proposition  \ref{t:psdup} provides  more detailed   NG sum upper bounds for $\thp$ and $\thpa$.

\begin{cor}\label{NG-Zp-prodth-both}
Let $G$ be a graph of order $n\ge 4$ such that both $G$ and $\ol G$ are connected.  Then:
\ben[$(1)$] 
\item\label{NGZpthbds}   $  n\le \thp(G)+ \thp(\ol G) \le 2n-o(n)$;\\
for $n\ne  5$, $ \max\{2n,3n-9\} \le \thp(G)\cdot \thp(\ol G) \le n^2-o(n^2)$\\
(for $n=  5$, $9 \le \thp(G)\cdot \thp(\ol G) \le  12$). 
\item\label{NGZpprodthabds} $  n-2\le \thpa(G)+ \thpa(\ol G) \le 2n-o(n)$ and \\$ \max\{n-1,2n-8\} \le \thpa(G)\cdot \thpa(\ol G) \le n^2-o(n^2)$. 

\item\label{NGZpprodthxbds} $  n\le \thpx(G)+ \thpx(\ol G)  \le 2n$ and for $n\ge 6$, $ 3n \le \thpx(G)\cdot \thpx(\ol G)  \le n^2$;\\  for $n=4$, $  \thpx(G)\cdot \thpx(\ol G) =9$  and for $n=5$, $ 12 \le \thpx(G)\cdot \thpx(\ol G)  \le 20$. 
\een
  For \eqref{NGZpthbds} and \eqref{NGZpprodthabds}, the NG upper bounds  are  tight and the NG  lower bounds  are  almost sharp. \\
 For \eqref{NGZpprodthxbds}, the NG upper bounds and the NG product lower bound   are   sharp and the NG sum lower   bound is  almost sharp. 
 \end{cor}

\bpf  The upper bounds are all immediate by Corollary \ref{NG-Zp-prodth-both-arb}  (the cases $n=5$ for $\thp(G)\cdot \thp(\ol G)$ and $n=4,5$ for $\thpx(G)\cdot \thpx(\ol G)$ are done in \cite{sage}). 
Since $c=\sqrt{|V(\fc_c)|}$,  Theorem \ref{p:Fc} shows the upper bounds in   \eqref{NGZpthbds} and \eqref{NGZpprodthabds} are tight and the upper bounds in \eqref{NGZpprodthxbds}  are sharp. 

The NG sum lower bounds were established in Corollary \ref{NG-Zp-prodth-both-arb}, and the NG product lower bounds were established in Lemma \ref{Lem:NG-Zp-prodth-lower-arb} because $\ol{K_{1,n-1}}$ is disconnected.   (See Remark \ref{rem:NG-Zp-prodth-lower-arb} for the small graph exceptions.)

Recall the graph $\ls_n=S(2,1,\dots,1)$ from Example \ref{ex:star+leaf}.  Since  $\Zp(\ls_n)=1$, $\ptp(\ls_n,1)=2$, and $\ptp(\ls_n,2)=1$, it follows that $\thpa(\ls_n)=2$ and $\thp(\ls_n)=\thpx(\ls_n)=3$.   
Furthermore, $\Zp(\ol{\ls_n})= n-3$,   $\ptp(\ol{\ls_n},n-3)=2$, and $\ptp(\ol{\ls_n},n-2)=1$, so $\thp(\ol{\ls_n})=n-1$, $\thpa(\ol{\ls_n})=n-2$,   and $\thpx(\ol{\ls_n})=n$. 
Thus $\thp(\ls_n)+\thp(\ol{\ls_n})=n+2$, $\thpa(\ls_n)+\thpa(\ol{\ls_n})=n$,  $\thpx(\ls_n)+\thpx(\ol{\ls_n})=n+3$, and all three NG sum lower bounds are almost sharp. 
 Moreover, $\thp(\ls_n)\cdot\thp(\ol{\ls_n})=3n-3$ and $\thpa(\ls_n)\cdot\thpa(\ol{\ls_n})=2n-4$, so the NG product lower bounds for $\thp$ and $\thpa$ are almost sharp.  Finally, $\thpx(\ls_n)\cdot\thpx(\ol{\ls_n})=3n$, so the NG product lower bound for $\thpx$ is sharp.  
 \epf

We end this section by describing the upper bounds $2n-o(n)$ for $\thp(G)+\thp(G)$ and $\thpa(G)+\thpa(G)$ in Corollary \ref{NG-Zp-prodth-both} more precisely.

\begin{prop}\label{t:psdup}
Let $G$ be a graph of order $ n\ge 4$ such that both $G$ and $\ol G$ are connected.  Then:
\ben[$(1)$] 
\item\label{NGZpthU}  $\thp(G)+\thp(\ol G)\le 2n-\log_4 n+2$. 
\item\label{NGZpprodthU}  $\thpa(G)+\thpa(\ol G)\le 2n-\log_4 n$. 
\een
\end{prop}

\bpf

By \cite[Proposition 2.9]{PSDthrottle},
$ \thp(G) \leq n - \alpha(G) + 1.$ 
Simply using this bound for both $G$ and $\ol G$, we get
\begin{equation}\label{eqn:thpsumbound}
\thp(G) + \thp(\ol G) \leq 2n - [\alpha(G) + \alpha(\ol G)] + 2.
\end{equation}
 The next NG bound on $\alpha$ is known \cite[Theorem 2.100]{NGsurvey13}:
\begin{equation}\label{Ramsey-bound}
\min\{a + b : R(a+1, b+1) > n\} \leq \alpha(G) + \alpha(\ol G),
\end{equation}
where $R(a+1, b+1)$ is the Ramsey number  for $a+1$ and $b+1$, i.e., the smallest integer $N$ such that every graph on $N$ vertices has an $(a+1)$-clique or a $(b+1)$-independent set.

Without loss of generality, suppose $a \geq b$. Due to a classic result of Erd\H{o}s and Szekeres \cite{ES47},
\[R(a+1, b+1) \leq R(a+1,a+1) \leq 4^a.\]
Therefore, the left-hand side of \eqref{Ramsey-bound} is at least
\[\min\{a + b : 4^a > n\} \geq \min\{a : 4^a > n\} \geq \log_4 n.\]
Combining this with  \eqref{eqn:thpsumbound} and \eqref{Ramsey-bound} gives the desired result:
\[\thp(G) + \thp(\ol G) \leq 2n - [\alpha(G) + \alpha(\ol G)] + 2 \leq 2n - \log_4 n + 2.\]

The second result follows in a similar manner from the fact that $\thpa(G) \leq n - \alpha(G)$ \cite[Proposition 2.77]{product2}.
\epf


\section{Power domination and Cops and Robbers} \label{s:PD+CR-NG-th}

In this section we establish Nordhaus-Gaddum bounds on throttling for two forcing parameters, namely power domination and Cops and Robbers. We combine the studies of these two forcing parameters because many bounds for both parameters are determined from domination bounds, including those in the next remark.  
For results that apply to both forcing parameters, we use $Y$ in place of the parameters $\pd$ or $c$.

\begin{rem}\label{r:thdom} 
    Because a dominating set is both a power dominating set and a capture set that has propagation time/capture time  at most one, we have the following upper bounds for every graph $G$.
\ben[$(1)$]
    \item $\thx(G)\le \gamma(G)+1$.
     \item $\thxa(G)\le \gamma(G)$.
      \item $\thxx(G)\le 2\gamma(G)$.
    \een
\end{rem}

  The next result presents some well-known  Nordhaus-Gaddam bounds on domination numbers.

\begin{thm}\label{NGdomthm}
  Let $G$ be a graph of order $n$. 
  \ben[$(1)$]
   \item\label{NGdomthm-d3}{\rm \cite[Theorem 2.1]{NGsurvey13}} Then $\gamma(G)+ \gamma(\Gc)\le n+1$ and $\gamma(G)\cdot \gamma(\Gc)\le n$. 
  \item\label{NGdomthm-NI} {\rm \cite[Theorem 2.9]{NGsurvey13}} If neither $G$ nor $\Gc$ has isolated vertices, then $\gamma(G)+\gamma(\Gc)\le \lf \frac n 2 \rf + 2$.  Furthermore, if $\gamma(G)+\gamma(\Gc)= \lf \frac n 2 \rf + 2$, then $\{\gamma(G),\gamma(\Gc)\} = \{\lf \frac n 2 \rf , 2\}$ or $G\cong K_3\cp K_3$. 
    \item\label{NGdomthm-d2}{\rm \cite[Theorem 2.10]{NGsurvey13}} If  $G$ and $\Gc$ are connected with $\delta(G),\delta(\Gc)\ge 2$ and $n\ge 10$, 
    then $\gamma(G)+\gamma(\Gc)\le \lf \frac {2n} 5\rf + 2$.
 \een
\end{thm}

\begin{rem}\label{r:conn-LB-pd-c} 
Let $G$ be a graph of order $n\ge 4$.   If  neither  $G$ nor $\Gc$ has a universal vertex, then  $\thx(G), \thxx(G), \thx(\Gc), \thxx(\Gc)\ge 3$ and $\thxa(G), \thxa(\Gc)\ge 2$.  

Now suppose  that $G$ has a universal vertex $u$.  Then $\thx(G)=\thxx(G)=2$, $\thxa(G)=1$, and $\ol G=K_1\du \ol{G-u}$. Since 
$Y(H)\geq2$ for any disconnected graph $H$, we have $ \thx(\ol G) \ge Y(\ol G)+1\ge 3$, $\thxa(\ol G) \ge Y(\ol G)\ge 2$, and   $ \thxx(\ol G)\ge \min\{2(1+1),n\}\ge 4$. 
\end{rem}

The examples below are used to establish sharpness of various NG 
bounds.

\begin{ex}\label{ex:star-pd-cr}
   Let 
   $n\ge 4$. Then $Y(K_{1,n-1})=1$, $\ptx(K_{1,n-1})=1$,  $Y(\ol{K_{1,n-1}})=2$, and $\ptx(\ol{K_{1,n-1}})=1$.  
   Thus \[\thx(K_{1,n-1})=2;\ \thxa(K_{1,n-1})=1;\ \thxx(K_{1,n-1})=2;\] 
   \[\thx(\ol{K_{1,n-1}})=3;\ \thxa(\ol{K_{1,n-1}})= 2;\ \thxx(\ol{K_{1,n-1}})=4.\]
\end{ex}

\begin{ex}\label{ex:Kn-pd-cr}   Let 
$n\ge 4$. Then $Y(K_n)=1$, $\ptx(K_n)=1$,  $Y(\ol{K_n})=n, $ and $\ptx(\ol{K_n})=0$.  Thus $\thx(K_n) =\thxx(K_n)=2$ and $\thx(\ol{K_n}) =\thxx(\ol{K_n})=n$.  
For $\thxa$, each graph must have an edge.  Since $Y(K_n-e)=1$, $\ptx( K_n-e)=1$, 
$Y(\ol{K_n-e})=n-1$, and $\ptx(\ol{K_n-e})=1$, we have $\thxa(K_n-e)=1$ and $\thxa(\ol{K_n-e})=n-1$.
\end{ex}

A vertex $u\in V(G)$ is a \emph{corner} if there is some  other
vertex $w$ 
such that $N[u]\subseteq N[w]$.  It is well known that any graph $G$ with $c(G)=1$ has a corner \cite[Lemma 2.1]{CRbook}.
\begin{ex}\label{ex:rK2-pd-cr}
   Let $r\ge 2$. Then $Y(rK_2)=r$, $\ptx(rK_2)=1$,  and $ \pd(\ol{rK_2})=1$. Since $\ol{rK_2}$ does not have a universal vertex,  $\ppt(\ol{rK_2})=2$. Since $\ol{rK_2}$ does not have a corner and any two  
   vertices dominate  $\ol{rK_2}$, we have $ c(\ol{rK_2})=2$  and $\capt(\ol{rK_2})=1$.  Thus 
\begin{center}
\begin{tabular}{ccc}
$\thx(rK_2)=r+1$;& \ $\thxa(rK_2)=r$;& $\thxx(rK_2)=2r$;\\
$\thx(\ol{rK_2})=3$;& $\thxa(\ol{rK_2})= 2$;& \ \ $\thpdx(\ol{rK_2})=3$, $\thcx(\ol{rK_2})=4$.
\end{tabular}
\end{center}
\end{ex}

\begin{ex}\label{r:pd-leaf} 
  Let $G$ be a  graph with a    leaf that does not have a universal vertex  or any isolated vertices.  Then $\pd(\Gc)=1$, $\ppt(\Gc)=2$,  
  $\thpd(\Gc)=3$, $\thpda(\Gc)=2$, and $\thpdx(\Gc)=3$. Furthermore, $c( \Gc)\le 2$ 
  and $\capt(\Gc,2)=1$, so $\thc(\Gc)=3$ and $\thca(\Gc)=2$.
\end{ex}

\begin{ex}\label{ex:HcircK1circK1-pd}
   Let $H$ be a connected graph, let $n=4|V(H)|$, and let $G=(H\circ K_1)\circ K_1$. Then $\thpda(G)=\frac n 2$ \cite{product1} and since Example~\ref{r:pd-leaf} applies to $G$,
   $\thpda(\Gc)=2$.
\end{ex}

 \begin{lem}\label{lem:thpd<=n+3 and 3n}
Let $G$ be a graph on $n\geq 4$ vertices  with $n$ even such that $\{\gamma(G),\gamma(\Gc)\}=\{\frac n 2,2\}$. Then $\thpdx(G)+\thpdx(\Gc)\le n+3$ and $\thpdx(G)\cdot\thpdx(\Gc)\le3n$.
\end{lem}
\bpf Observe that neither $G$ nor $\Gc$ has a universal vertex or an isolated vertex.
Without loss of generality, let $\gamma(G)=\frac n 2$.  Recall that $\thpdx(G)\le n$.   By \cite[Theorem 4.24]{core-graph-dom-book}, every component of $G$ must be a $4$-cycle or a corona $H\circ K_1$ for some graph $H$. If any component of $G$ is of the form $H\circ K_1$, then $\thpdx(\Gc)=3$ by Example \ref{r:pd-leaf}. Hence, $\thpdx(G)+\thpdx(\Gc)\le n+3$ and $\thpdx(G)\cdot\thpdx(\Gc)\le 3n$. Next, consider the case $G=rC_4$ (the disjoint union of $r$ copies of $C_4$). 
Observe that  $\thpdx(rC_4)=3r=\frac {3n}4$ and $\thpdx(\ol{rC_4})\le 4$ because  two vertices from distinct copies of $C_4$ in $rC_4$ form a dominating set for $\ol{rC_4}$; thus $\thpdx(G)+\thpdx(\Gc)\le \frac{3n}4+4\le n+3$ and $\thpdx(G)\cdot\thpdx(\Gc)\le\frac{3n}{4}\cdot4=3n$.
\epf

\begin{thm}\label{NG-pd-cr}  Let $G$ be a graph of order $n\ge 5$. 

For power domination number:
    \ben[$(1)$]

    \item $5\le \thpd(G)+ \thpd(\ol G) \le n+2$ and $6\le \thpd(G)\cdot \thpd(\ol G)\le 2n$.

    \item $3\le \thpda(G)+ \thpda(\ol G) \le n$ and $2\le \thpda(G)\cdot \thpda(\ol G) \le n$. 

    \item\label{it:thpd} $6\le \thpdx(G)+ \thpdx(\ol G) \le n +3$   
    and $8\le \thpdx(G)\cdot \thpdx(\ol G)\le    4n$. 

For cop number:
            \item $5\le \thc(G)+ \thc(\ol G) \le n+2$ and $6\le \thc(G)\cdot \thc(\ol G)\le 2n$.

    \item $3\le \thca(G)+ \thca(\ol G) \le  n$ and $2\le \thca(G)\cdot \thca(\ol G) \le  n$. 

    \item $6\le \thcx(G)+ \thcx(\ol G) \le n +4$  and $8\le \thcx(G)\cdot \thcx(\ol G)\le     4n$. 
    \een
 All these bounds are sharp except $\thpdx(G)\cdot \thpdx(\ol G)\le   4n$,  which is sharp only for $K_3\cp K_3$ (and may not be tight  as $n\to\infty$).   

If  neither $G$ nor $\Gc$ has an isolated vertex, then \beq\label{eq:noisol} \thx(G)+\thx(\Gc)\le \lf \frac{n}2\rf+4 \ \mbox{ and } \ \thx(G)\cdot\thx(\Gc)\le n+ \lf \frac n 2\rf +3.\eeq
\end{thm}
\bpf     
 Remark \ref{r:conn-LB-pd-c} establishes the lower bounds and Example \ref{ex:star-pd-cr} shows that they are sharp.

Upper bounds for NG products:  Remark \ref{r:thdom} and Theorem \ref{NGdomthm}\eqref{NGdomthm-d3} imply $\thxa(G)\cdot\thxa(\Gc)\le n$ and $\thxx(G)\cdot\thxx(\Gc)\le 4n$.  If $G$ has an isolated vertex, then $\gamma(\Gc)=1$ so $ \thx(\Gc)=2$ and $\thx(G)\cdot\thx(\Gc)\le 2n$.  If neither $G$ nor $\Gc$ has an isolated vertex, then by Remark \ref{r:thdom} and Theorem \ref{NGdomthm}\eqref{NGdomthm-NI}, 
\[\thx(G)\cdot\thx(\Gc)\le(\gamma(G)+1)(\gamma(\Gc)+1)=\gamma(G)\cdot\gamma(\Gc)+\gamma(G)+\gamma(\Gc)+1\le n+ \lf \frac n 2\rf +3 \le 2n \] 
where the last inequality is true because  $n\ge 5$.  This also establishes the second statement in \eqref{eq:noisol}.

Upper bounds for NG sums: Whenever neither $G$ nor $\Gc$ has an isolated vertex, Remark \ref{r:thdom} and Theorem \ref{NGdomthm}\eqref{NGdomthm-NI}  imply $\thx(G)+\thx(\Gc)\le \lf \frac{n}2\rf+4$ (which establishes the first statement in \eqref{eq:noisol}),  $\thxa(G)+\thxa(\Gc)\le \lf \frac{n}2\rf+2$, and $\thxx(G)+\thxx(\Gc)\le n+4$. Now suppose $G$ has an isolated vertex $u$, so $u$ is a universal vertex in $\Gc$. 
  In this case, $\thx(G)+\thx(\Gc)\le n+2$, $\thxa(G)+\thxa(\Gc)\le n$, and $\thxx(G)+\thxx(\Gc)\le n+2$ since $\thx(\Gc)= \thxx(\Gc)=2$,  $\thxa(\Gc)=1$, $\thx(G)\le n$, $\thxx(G)\le n$,  and $\thxa(G)\le n-1$  by Remark~\ref{NG-u}. 
  The larger of the two bounds is selected for each type of throttling.

Thus every bound except the NG sum upper bound for $\thpdx$ has been established.  To establish $n+3$ as the bound, we analyze the proof of the $n+4$ NG sum upper bound for $\thpdx$.  Recall that if $G$ or $\Gc$ has an isolated vertex, then $\thxx(G)+\thxx(\Gc)\le n+2$, so assume neither $G$ not $\Gc$ has an isolated vertex.  Then by Theorem \ref{NGdomthm}\eqref{NGdomthm-NI}, $\thpdx(G)+\thpdx(\Gc)\le 2(\lf \frac{n}2\rf+2)$.  Note first that if $n$ is odd then  $ 2(\lf \frac{n}2\rf+2)\le n+3$, so assume $n$ is even. If $\gamma(G)+\gamma(\Gc)< \frac n 2  + 2$, then $\thpdx(G)+\thpdx(\Gc)< 2( \frac{n}2+2)=n+4$. So  assume $\gamma(G)+\gamma(\Gc)= \frac n 2  + 2$.  Then by Theorem \ref{NGdomthm}\eqref{NGdomthm-NI},
$\{\gamma(G),\gamma(\Gc)\} = \{ \frac n 2  , 2\}$ (since $n\ne 9$ because $n$ is even), and we obtain the  stated upper bound for $\thpdx(G)+\thpdx(\Gc)$ by Lemma~\ref{lem:thpd<=n+3 and 3n}.  

A similar approach establishes that the bound $\thpdx(G)\cdot \thpdx(\ol G)\le 4n$ is not sharp except for $K_3\cp K_3$. From Remark~\ref{r:thdom} and Theorem~\ref{NGdomthm}\eqref{NGdomthm-d3}, we see that $\thpdx(G)\cdot \thpdx(\ol G)=4n$ only if $\gamma(G)\cdot\gamma(\ol G)=n$. By \cite[Theorem 16.7]{core-graph-dom-book}, either $\{\gamma(G), \gamma(\ol G)\} \in \{\{n,1\}, \{\frac{n}{2}, 2\}\}$ or $G=K_3\Box K_3$. However, if $\{\gamma(G), \gamma(\ol G)\}=\{n,1\}$, then $\thpdx(G)\cdot \thpdx(\ol G)\leq2n$, and if $\{\gamma(G), \gamma(\ol G)\}=\{\frac{n}{2}, 2\}\}$, then $\thpdx(G)\cdot \thpdx(\ol G)\leq3n$ by Lemma~\ref{lem:thpd<=n+3 and 3n}. Thus $\thpdx(G)\cdot \thpdx(\ol G)<4n$ except for $K_3\cp K_3$, and $\thpdx(K_3\cp K_3)\cdot \thpdx(\ol{K_3\cp K_3})= 6\cdot 6=4\cdot 9$ \cite{sage}.

  Example \ref{ex:Kn-pd-cr} shows  both of the NG upper bounds for $\thx$  are sharp and the NG sum upper bound for  $ \thxa(G)$ is sharp. 
 Example \ref{ex:rK2-pd-cr} shows that  the NG product upper bound for $\thxa$ is sharp, both NG 
 upper bounds  for $\thcx$ are sharp,  and the NG sum upper bound for $\thpdx$ is sharp (note the order there is $2r$). 
\epf

\begin{rem}
    Theorem \ref{NG-pd-cr} is stated for graphs of order $n\ge 5 $.  However, all  bounds apply to graphs of order $n=4$  except the NG product upper bound on $\thx$. The only pairs of graphs $G,\Gc$ of order $n=4$ that violate   $\thx(G)\cdot \thx(\Gc)\le 2\cdot 4$ are $\thx(P_4)\cdot \thx(P_4)=9$ and $\thx(C_4)\cdot \thx(2K_2)=9$ \cite{sage}. 
\end{rem}

 The next result shows that any upper bound on $\thpdx(G)\cdot \thpdx(\ol G)$ for graphs $G$ of order $n$  is at least $\frac{24n}{7}\approx 3.42857n$.

 \begin{cor}\label{NG-pd-cr-theta} 
   For $r\ge 1$ there exists a graph $G$ of order $n=7r$ such that $\thpdx(G)\cdot \thpdx(\ol G)=24r=\frac{24n}{7}$.  Thus  the NG product upper bound for $\thpdx$ is $\Theta(n)$. 
\end{cor}
\bpf
    The NG product upper bound for $\thpdx$ is  less than $4n$  by    Theorem \ref{NG-pd-cr}, except for $G = K_3 \cp K_3$. 
    
    Let $G=rC_3\sqcup rC_4$ be a graph on $n=7r$ vertices. Then $\gamma(G)=3r$, $\pd(G)=2r$, and $\ppt(G)=2$. Hence, $\thpdx(G)=\min(3r(1+1),2r(2+1))=6r$. The complement $\Gc$ is the join of $r$ copies of $3K_1$ and $r$ copies of $2K_2$, so $\pd(\Gc)=2=\gamma(\Gc)$. 
Thus $\thpdx(G)\cdot\thpdx(\ol G)=6r\cdot4=\frac{24n}{7}$. \epf

Next we establish Nordhaus-Gaddum bounds for power domination  and Cops and Robbers throttling when both graphs are required to be connected.

 \begin{ex}\label{ex:thpd}
    Consider the graph $G_0$ of order $n$ obtained by adding a new leaf to one leaf of $H\circ2K_1$, where $H$ is a connected graph of order at least $4$. It was shown in \cite[Example 11.57]{HLSbook} that $\thpd(G_0)=\lf\frac{n}{3}\rf+2$. Moreover, $\gamma(\ol{G_0})=2$ and $\thpd(\ol{G_0})=3$, so $\thpd(G_0)+ \thpd(\ol{G_0}) =\lf\frac n 3\rf +5$. Also $\thpd(G_0)\cdot \thpd(\ol{G_0})=3\lp\lf\frac{n}{3}\rf+2\rp= n+o(n)$.
\end{ex}

\begin{thm}\label{t:NG-pd-c-U} Let $G$ be a graph of order $n$ such that both $G$ and $\ol G$ are connected.  

If  $n\ge 4$ and $n\ne 9,10,11$, then  \[\thpd(G)+ \thpd(\ol G) \le  \lf\frac n 3\rf +5\] and this bound is  sharp. Furthermore, $\thpd(G)+ \thpd(\ol G) \le  \lf\frac n 3\rf +6$ for $n=9,10$ and $\thpd(G)+ \thpd(\ol G) \le  \lf\frac n 3\rf +7$ for $n=11$.

For all $n\geq4$,  \[\thpd(G) \cdot \thpd(\ol G) \le  n+2(1+\ln n)\sqrt{n}+4.\]
\end{thm}
\bpf
Since $G$ and $\ol G$ are connected, we have $\gamma(G)\geq2$ and $\gamma(\ol G)\geq2$.  Suppose $\gamma(\ol G)=2$; this implies $\thpd(\ol G) =3$.
It was shown in \cite
{ourpaper1} that $\thpd(G)\leq\lf\frac{n}{3}\rf+2$, so $\thpd(G)+ \thpd(\ol G) \le  \lf\frac n 3\rf +5$ and $\thpd(G) \cdot \thpd(\ol G) \le  n+6$. Similarly, we obtain the same upper bounds if $\gamma(G)=2$. 
Hence, it remains to consider  the case where $\gamma(G)\geq3$ and $\gamma(\ol G)\geq3$.
Under this assumption, we have $\gamma(G)+\gamma(\ol G)\leq\min\{\delta(G),\delta(\ol G)\}+3$ \cite[Theorem~2.7]{NGsurvey13}.

  If $\min\{\delta(G),\delta(\ol G)\}\leq\lf\frac{n}{3}\rf$, then $\thpd(G)+ \thpd(\ol G)\leq\gamma(G)+\gamma(\ol G)+2\leq\lf\frac{n}{3}\rf+5$. Hence, in the remainder of the proof  of the NG sum upper bound, we assume that $\delta(G),\delta(\ol G)>\lf\frac{n}{3}\rf$. 
Since $\gamma(G)\leq\left(\frac{1+\ln(\delta(G)+1)}{\delta(G)+1}\right)n$ for every graph $G$ with no isolated vertices \cite[Theorem~16]{henning}, we have $\max\{\gamma(G),\gamma(\ol G)\}\le\left(\frac{1+\ln(\lf n/3\rf+{2})}{\lf n/3\rf+{2}}\right)n$, where the  inequality follows from the fact that $f(x) = \frac{1 + \ln(x)}{x}$ is a decreasing function on $x\geq 1$ and $\delta(G)\ge\lf\frac n 3 \rf+1$. As a result, \[\thpd(G)+ \thpd(\ol G)\leq 2\left(\frac{1+\ln(\lf n/3\rf+{2})}{\lf n/3\rf+2}\right)n +2,\]
 which is numerically verified to be less than $\lf n/3\rf+5$ for all $n\geq 63$.

Now assume $15\le n\le 62$. Since $n\ge 15$,  $\delta(G),\delta(\Gc)\ge 6$.  By \cite[Theorem~2.16]{NGsurvey13}, $\gamma(G)+\gamma(\Gc)\le \lf\frac{6n}{17}\rf+2$, which implies 
\[\thpd(G)+ \thpd(\ol G)\leq \lf\frac{6n}{17}\rf+4\le \lf\frac n 3\rf+5\]
because $15\le n\le 62$.  The case $n=12,13,14$ is similar using $\gamma(G)+\gamma(\Gc)\le \lf\frac{5n}{14}\rf+2$ when $\delta(G),\delta(\Gc)\ge 5$ and $n\ne 16$ by \cite[Theorem~2.15]{NGsurvey13}.

Finally we consider $4\le n\le 11$.  For all pairs $G$ and $\Gc$ where both are connected, $\thpd(G)+ \thpd(\ol G)\leq  \lf\frac n 3\rf+5$ for $n=4,5,6,7,8$ is verified in \cite{sage}. Since $\min\{\delta(G),\delta(\ol G)\}+\max\{\Delta(G),\Delta(\ol G)\}=n-1$, $\min\{\delta(G),\delta(\ol G)\}\le \lf\frac{n-1}2\rf$.  Thus $\thpd(G)+ \thpd(\ol G) \le  \lf\frac {n-1} 2\rf +5$.  Computation shows that $\lf\frac {n-1} 2\rf +5\le \lf\frac {n} 3\rf +5+x$ with $x=1$ for $ n=9,10$ and $x=2$ for $n=11$.

This completes the proof of the NG sum upper bound for $\thpd$, and Example~\ref{ex:thpd} shows  that the bound for $n\geq 12$ is sharp.


As for the NG product upper bound, \[\thpd(G)\cdot\thpd(\ol G)\leq(\gamma(G)+1)(\gamma(\ol G)+1)\leq n+\gamma(G)+\gamma(\ol G)+1\leq n+\min\{\delta(G),\delta(\ol G)\}+4\] 
 since $\gamma(G)+\gamma(\ol G)\leq\min\{\delta(G),\delta(\ol G)\}+3$.
If $\min\{\delta(G),\delta(\ol G)\}\leq\sqrt{n}$ (which  is always true if $n\le 6$), then our result is established.
So assume $\min\{\delta(G),\delta(\ol G)\}>\sqrt{n}$. Then $\max\{\gamma(G),\gamma(\ol G)\}<\left(\frac{1+\ln(\sqrt{n}+1)}{\sqrt{n}+1}\right)n<(1+\ln n)\sqrt{n}$ since $n\ge 7$, which implies that $\thpd(G)\cdot\thpd(\ol G)<n+2(1+\ln n)\sqrt{n}+4$. 
\epf

The following theorem supplies an upper bound for the NG sum of $\thc$ but is likely not asymptotically tight.  Note that the largest known example has $\thc(G)=\Omega(n^{2/3})$ \cite{cop-throttle2}. 

\begin{thm}\label{cr-thc-ubd} {\rm\cite{cop-throttle2}}
    If $G$ is a connected graph of order $n$, then $ \thc(G) \le \frac{n}{(\log n)^{\frac{1}{2}-o(1)}}$.
\end{thm}
\begin{ex}\label{ex:leaf-star-pd-cr}
Let $\ls_n$ be the graph defined in Example~\ref{ex:star+leaf}.  It is straightforward to verify that $c(\ls_n)=c(\ol{\ls_n})=1$ and $\capt(\ls_n)=\capt(\ol{\ls_n})=2$, so $\thcx(\ls_n)=\thcx(\ol{\ls_n})=3$. Furthermore, observe that each of $\ls_n$ and $\ol{\ls_n}$ has a leaf, so Example~\ref{r:pd-leaf} justifies  the remaining values in the next statement:
\[\thx(\ls_n)=\thx(\ol{\ls_n})=3,\ \thxa(\ls_n)=\thxa(\ol{\ls_n})=2,\ \thxx(\ls_n)=\thxx(\ol{\ls_n})=3.\]
\end{ex}

\begin{thm}\label{NG-pd-cr-con} Let $G$ be a graph of order $n\ge 5$  such that both $G$ and $\ol G$ are connected.
  For power domination number:
    \ben[$(1)$]
        \item\label{thpd-conn} $6\le \thpd(G)+ \thpd(\ol G) \le   \lf\frac n 3\rf +5$ $(n\ge  12)$  
        and $9\le \thpd(G)\cdot \thpd(\ol G)\le n+o(n)$. 
        The exact bounds are sharp and the asymptotic bound is tight.
    \item\label{thpda-conn} $4\le \thpda(G)+ \thpda(\ol G) \le \lf\frac n 2\rf +2$ and $4\le \thpda(G)\cdot \thpda(\ol G)\le n$.
    All bounds are sharp.   
    \item\label{thpdx-conn} $6\le \thpdx(G)+ \thpdx(\ol G) \le \frac {6n}7  +3$ $ (n\ge 16)$ and $9\le \thpdx(G)\cdot  \thpdx(\ol G)\le 4n$.  Both lower bounds and the NG sum upper bound are sharp.

   \hspace{-10mm} For cop  number :
   
        \item\label{thc-conn} $6\le \thc(G)+ \thc(\ol G) \le \min\lp\frac{2n}{(\log n)^{\frac{1}{2}-o(1)}},\lfloor \frac{n}{2} \rfloor + 4\rp$  and $9\le \thc(G)\cdot \thc(\ol G)\le n + \lfloor \frac{n}{2} \rfloor + 3$. Both lower bounds are sharp.

    \item\label{thca-conn} $4\le \thca(G)+ \thca(\ol G) \le \lf\frac n 2\rf +2$ and $4\le \thca(G)\cdot \thca(\ol G) \le n$.  Both lower  bounds are sharp. 

    \item\label{thcx-conn} $6\le \thcx(G)+ \thcx(\ol G) \le n+4$  and $9\le \thcx(G)\cdot \thcx(\ol G)\le  4n   $.  Both lower bounds are sharp.

    \een

\end{thm}
\bpf   

Lower bounds:  When both $G$ and $\Gc$ are connected, neither has a universal vertex. Thus Remark \ref{r:conn-LB-pd-c} establishes the lower bounds.   Example \ref{ex:leaf-star-pd-cr} shows that all lower bounds are sharp.

For the upper bounds in \eqref{thpd-conn} and \eqref{thc-conn}: Theorem \ref{t:NG-pd-c-U} establishes the upper bounds for $\thpd$ and  
 Example \ref{ex:thpd} shows those  upper bounds are tight. Since neither $G$ nor $\ol G$ have isolated vertices, \eqref{eq:noisol} and Theorem \ref{NG-pd-cr} 
establish the upper bounds for $\thc$.

  For the upper bounds in \eqref{thpda-conn} and \eqref{thca-conn}:
  Parts \eqref{NGdomthm-NI} and \eqref{NGdomthm-d3} of Theorem \ref{NGdomthm}
  show that $ \thxa(G)+ \thxa(\ol G) \le\gamma(G)+\gamma(\ol G)\le \lf\frac n 2\rf +2$ and $ \thxa(G)\cdot \thxa(\ol G) \le\gamma(G)\cdot\gamma(\ol G)\le n$, 
  respectively. 
  Example \ref{ex:HcircK1circK1-pd} shows that 
 the upper bounds for  $\thpda$ are sharp. 

Now we turn our attention to the upper bounds in \eqref{thpdx-conn} and \eqref{thcx-conn}.
   For $n\ge 10$, Theorem \ref{NGdomthm}\eqref{NGdomthm-d2}  shows $\thxx(G)+\thxx(\Gc)\le 2\lp\lf \frac{2n}5\rf+2\rp$ whenever $\delta(G),\delta(\Gc)\ge 2$. For $n\ge 16$, we have $2\lp\lf \frac{2n}5\rf+2\rp\le \frac {6n}7  +3$.

 Now suppose (without loss of generality) that $\delta(G)=1$.  Then $\Delta(\Gc)=n-2$, which implies $\pd(\Gc)=1$, $\ppt(\Gc)=2$, 
 and $\thpdx(\Gc)=3$.  Since $G$ is connected, $\thpdx(G)\le \frac{6n}7$ \cite{ourpaper1}.  Therefore, $\thpdx(G)+\thpdx(\Gc)\le \frac{6n}7+3$ regardless of the value of $\delta(G)$ (provided $n\ge 16)$.  The graph $G$ in Example 2.65 in \cite{product2} is a connected graph of order $n$ with $\thpdx(G)=\frac {6n}7$.  Since $G$ has a leaf, $\thpdx(\Gc)=3$.  Thus the NG sum upper bound for  $\thpdx$ is sharp.  
 
The NG product upper bound for $\thpdx$ and the upper bounds for $\thcx$ follow from the corresponding bounds in Theorem \ref{NG-pd-cr}.
\epf

 In the proof of Theorem \ref{NG-pd-cr-con}, it was noted that the graph $G$ in \cite[Example 2.65]{product2} is a connected graph of order $n$ with $\thpdx(G)=\frac {6n}7$ and $\thpdx(\Gc)=3$.  Thus  any upper bound on $\thpdx(G)\cdot \thpdx(\ol G)$ for graphs $G$ of order $n$ with both $G$ and $\ol G$ connected is at least $\frac{18n}{7}\approx 2.5714n$.

 \begin{cor}\label{NG-pd-theta-con}
   The NG product upper bound for $\thpdx$ is $\Theta(n)$. 
\end{cor}

\section*{Acknowledgements}
This research began and was subsequently continued at the American Institute of Mathematics with support from NSF DMS grant  2015462.  The authors thank AIM and NSF. Ryan Blair was supported in part by NSF grant DMS-2424734. Veronika Furst was supported in part by NSF grant DMS-2331072.



\begin{thebibliography}{99}

\bibitem{AIM08} AIM Minimum Rank -- Special Graphs Work Group  
{\scriptsize(F.~Barioli, W.~Barrett,  S.~Butler,    S.M.~Cioaba, D.~Cvetkovi\'c,   S.M.~Fallat, C.~Godsil,  
W.~Haemers,   L.~Hogben,  R.~Mikkelson,  S.~Narayan,  O.~Pryporova,   
I.~Sciriha,  W.~So,   D.~Stevanovi\'c,  H.~van der Holst, K.~Vander Meulen, and A.~Wangsness Wehe}).  
Zero forcing sets and the minimum rank   of graphs.   {\em Linear Algebra Appl.}  428 (2008), 1628--1648, .

\bibitem{product2} S.E.~Anderson, K.L.~Collins, D.~Ferrero, L.~Hogben,  C.~Mayer, A.N.~Trenk, and S.~Walker.  Product throttling.  In {\em Research Trends in Graph Theory and Applications} (D.~Ferrero, L.~Hogben, S.R.~Kingan, and G.L.~Matthews, editors), Association for Women in Mathematics Series, Springer, New York, 2021, pp. 11--50. 

 \bibitem{product1} S.E.~Anderson, K.L.~Collins, D.~Ferrero, L.~Hogben,  C.~Mayer, A.N.~Trenk, and S.~Walker.  Product throttling for power domination.   \emph{Australas. J. Combin.}  {85} (2023), 248--272. 

\bibitem{NGsurvey13}	M. Aouchiche and P. Hansen. A survey of Nordhaus--Gaddum type relations.  
\emph{Discrete Appl. Math.} 161 (2013), 466--546.

\bibitem{smallparam} 	F.~Barioli, W.~Barrett, S.~Fallat, H.T.~Hall, 
L.~Hogben, B.~Shader, P.~van den Driessche, and H.~van der Holst.  
Zero forcing parameters and minimum rank problems.  
{\em Linear Algebra Appl.}, 433 (2010), 401--411.

\bibitem{PD-NG-REUF} K.F. Benson, D. Ferrero, M. Flagg, V. Furst, L. Hogben, and V. Vasilevska. Note on Nordhaus-Gaddum problems for power domination.  {\it Discrete Appl. Math.},  {251} (2018),  103--113.

\bibitem{ourpaper1} R.~Blair, G.~Elvin, V.~Furst, L.~Hogben, N.~Sahajpal, and T.W.H.~Wong.  Sharp bounds for product and sum throttling numbers.  Under review.  \url{https://arxiv.org/abs/2501.03472}.

\bibitem{BC79} B. Bollob\'as and E. J. Cockayne. Graph-theoretic parameters concerning domination, independence, and irredundance. \emph{J. Graph Theory}, {\bf 3} (1979), 241-249.

\bibitem{cop-throttle2} A.~Bonato, J.~Breen, B.~Brimkov, J.~Carlson, S.~English, J.~Geneson, L.~Hogben, K.E.~Perry, and C.~Reinhart. Optimizing the trade-off between number of cops and capture time in Cops and Robbers.  {\em J.~Comb.} 13 (2022), 79--103.  

\bibitem{CRbook} A.~ Bonato and R.J.~Nowakowski. \emph{The game of Cops and Robbers on graphs.} American
Mathematical Society, Providence, 2011.

 

\bibitem{BY13} S.~Butler and M.~Young. Throttling zero forcing propagation speed on graphs. {\em Australas.~J.~Combin.} 57 (2013), 65--71.


\bibitem{powerdom-throttle}B.~Brimkov, J.~Carlson, I.V.~Hicks, R.~Patel, and L.~Smith. Power domination throttling.  
{\em Theoret.~Comput.~Sci.} 795 (2019), 142--153.


\bibitem{PSDthrottle}
J.~Carlson, L.~Hogben, J.~Kritschgau, K.~Lorenzen, M.S.~Ross, S.~Selken, V.~Valle Martinez. Throttling positive semidefinite zero forcing propagation time on graphs.
\emph{Discrete Appl. Math} 254 (2019), 33--46.

\bibitem{CH74} E. J. ~Cockayne, S. T. ~Hedetniemi. Independence graphs.In {\em Proceedings of the {F}ifth {S}outheastern {C}onference on
        {C}ombinatorics, {G}raph {T}heory and {C}omputing ({F}lorida
              {A}tlantic {U}niv., {B}oca {R}aton, {F}la., 1974)}, \emph{Congress. Numer.} X (1974),   471–491.

\bibitem{C22} E.~Conrad.  Positive Semidefinite Initial Cost Product Throttling. Available at \url{https://arxiv.org/abs/2207.02795}.

\bibitem{DVV16} P.~Dorbec,  S.~Varghese, A.~Vijayakumar. Heredity for generalized power domination. {\em Discrete Math.~Theor.~Comput.~Sci.}, 18 (2016), Paper No. 5, 11 pp.


\bibitem{ES47} P. Erd\H{o}s. Some remarks on the theory of graphs. {\em Bull. Amer. Math. Soc.} {53} (1947), 292--294. 

\bibitem{core-graph-dom-book} T.W.~Haynes, S.T.~Hedetniemi, M.A.~Henning. \emph{Domination in graphs—core concepts}.
Springer, Cham, 2023.

\bibitem{henning} M.A.~Henning. Bounds on domination parameters in graphs: a brief survey.
\emph{Discuss. Math. Graph Theory} 42 (2022),  665--708.

\bibitem{sage} L.~Hogben. {\em Sage} code verifications for this paper. Available at  \url{https://aimath.org/~hogben/NGthrottling.ipynb} (code) and \url{https://aimath.org/~hogben/NGthrottling_Sage.pdf} (PDF).

\bibitem{HLSbook} L. Hogben, J.C.-H. Lin, B.L. Shader.  {\em Inverse Problems and Zero Forcing for Graphs}.   {Mathematical Surveys and Monographs} {270}, American Mathematical Society, Providence, RI,  2022. 

\bibitem{Ore}
0.~Ore. \emph{Theory of Graphs}. Amer. Math. Soc. Colloq. Publ. XXXVIII, American Mathematical Society, Providence, RI, 1962.

\bibitem{SMH19} L.A.~Smith, D.J.~Mikesell, I.V.~Hicks.  An integer program for positive semidefinite zero forcing in graphs.
\emph{Networks}, 76 (2020), 366--380.
\bibitem{Zhao} M. Zhao, L. Kang, G.J. Chang.  Power domination in graphs.  \emph{Discrete Math.} {306} (2006), 1812--1816.


\end{thebibliography}
\end{document}